# Calculation of oblate spheroidal wave functions with complex argument

Arnie L. Van Buren

Cary, North Carolina

August 31, 2020

MSC -class: 3310

## Abstract

A previous article showed that alternative expressions for calculating oblate spheroidal radial functions of both kinds $R_{ml}^{(1)}(-ic,i\xi)$ and $R_{ml}^{(2)}(-ic,i\xi)$ can provide accurate values over very large parameter ranges using double precision arithmetic, even where the traditional expressions fail. The size parameter $c$ was assumed real. This paper considers the case where $c = c_r + ic_i$ is complex with an imaginary part $c_i$ often used to represent losses in wave behavior. The methods for $c$ real modified to complex arithmetic work reasonably well as long as $c_i$ is very small. This paper describes the substantial changes necessary to obtain useful results for larger values of $c_i$. It shows that accurate eigenvalues can usually be obtained even though the matrix methods used to obtain them for $c$ real provide increasingly inaccurate values, primarily for those with relatively small magnitude, as $c_i$ increases. It also shows that some of the eigenvalues can be prolate-like with values that are well approximated using asymptotic estimates for prolate eigenvalues where $c$ is replaced with $-ic$. A method to order the eigenvalues is presented. The modifications necessary to compute accurately the radial and angular functions for complex $c$ are discussed. A resulting Fortran program coblfcn provides useful function values for a reasonably wide range of $c$, $m$ and $\xi$ when using double precision arithmetic. The results can be improved by using quadruple precision for the Bouwkamp procedure to ensure accurate double precision eigenvalues. Further improvement is obtained using full quadruple precision. Coblfcn is freely available at www.mathieuandspheroidalwavefunctions.com

## 1  Introduction

The scalar Helmholtz wave equation for steady waves, $(\nabla^2 + k^2)\Psi = 0$, where $k = 2\pi/\lambda$ and $\lambda$ is the wavelength, is separable in the oblate spheroidal coordinates $(\xi, \eta, \varphi)$, with $0 \leq \xi \leq \infty$, $-1 \leq \eta \leq 1$, and $0 \leq \varphi \leq 2\pi$. The factored solution is $\Psi_{ml}(\xi, \eta, \varphi) = R_{ml}(-ic, i\xi) S_{ml}(-ic, \eta) \Phi_m(\varphi)$, where $R_{ml}(-ic, i\xi)$ is the radial function, $S_{ml}(-ic, \eta)$ is the angular function, and $\Phi_m(\varphi)$ is the azimuthal function. Here $c = ka/2$, where $a$ is the interfocal distance of the elliptic cross section of the spheroid. The radial function of the first kind $R_{ml}^{(1)}(-ic,i\xi)$ and the radial function of the second kind $R_{ml}^{(2)}(-ic,i\xi)$ are the two independent solutions to the second order radial differential equation resulting from the separation of variables. These solutions are dependent on four parameters $(m,l,c,\xi)$ and an eigenvalue (separation constant) $A_{ml}(-ic)$. Similarly, $S_{ml}^{(1)}(-ic,\eta)$ and $S_{ml}^{(2)}(-ic,\eta)$ are the two independent solutions to the second order angular differential equation resulting from the separation of variables. In the following discussion the order *m* is either zero or a positive integer with the degree *l* equal to *m, m+1, m+2, ….*

Oblate spheroidal functions are used in solving boundary value problems of radiation, scattering, and propagation of scalar and vector acoustic and electromagnetic waves in oblate spheroidal coordinates.

A previous paper [1] shows that alternative expressions for calculating the oblate spheroidal radial functions of both kinds $R_{ml}^{(1)}(-ic,i\xi)$ and $R_{ml}^{(2)}(-ic,i\xi)$ can provide accurate values over very large parameter ranges using double precision arithmetic, even where the traditional expressions fail. It describes some of the features of a new Fortran computer program oblfcn that calculates the oblate spheroidal angular and radial functions using a combination of both the traditional and the alternative expressions. Here *c* was assumed real. In this paper $c = c_r + ic_i$ is assumed complex with the imaginary component $c_i$ arising from the imaginary component of *k* that often represents losses in wave behavior. It is reasonably straightforward to convert oblfcn to complex arithmetic. The resulting program provides accurate oblate function values as long as $c_i$ is less than about 5. Making modifications to accommodate larger values of $c_i$ is much more involved. This paper discusses those modifications and the resulting Fortran program coblfcn. It concludes with a summary.

## 2  Angular functions of the first kind

The oblate angular function of the first kind $S_{ml}^{(1)}(-ic,\eta)$ is expressed [see for example ref. 2, p. 16] in terms of the corresponding associated Legendre functions of the first kind by

$$S_{ml}^{(1)}(-ic,\eta) = \sum_{n=0,1}^{\infty} {'} d_n(-ic|ml) P_{m+n}^{m}(\eta), \qquad (1)$$

where the prime sign on the summation indicates that *n* = 0, 2, 4,... if *l − m* is even or *n* = 1, 3, 5, ... if *l − m* is odd. The following three term recursion formula relates successive expansion coefficients $d_{n-2}, d_n,$ and $d_{n+2}$ for given values of *l, m,* and *c*:

$$\alpha_n d_{n+2} + (\beta_n - \lambda_{ml}) d_n + \gamma_n d_{n-2} = 0, \qquad (2)$$

where

$$\alpha_n = -\frac{(2m+n+2)(2m+n+1)c^2}{(2n+2m+3)(2n+2m+5)},$$

$$\beta_n = \left[(n+m)(n+m+1) - \frac{2(n+m)(n+m+1) - 2m^2 - 1}{(2n+2m+3)(2n+2m-1)} c^2\right],$$

$$\gamma_n = -\frac{n(n-1)c^2}{(2n+2m-3)(2n+2m-1)}. \qquad (3)$$

Use of this formula to calculate the expansion coefficients requires a value for the separation constant or eigenvalue $\lambda_{ml}(-ic)$, which is chosen to ensure nontrivial convergent solutions for $S_{ml}^{(1)}(-ic,\eta)$.

Oblfcn eigenvalues are more difficult to obtain that prolate eigenvalues. In the prolate case, one uses traditional approximations [3] for the lowest eigenvalues beginning with $l = m$. One then uses extrapolations of previous eigenvalues to obtain approximations for higher eigenvalues. The Bouwkamp procedure [4] refines the approximations to provide accurate eigenvalues. For the oblate case this approach does not work well unless $c$ is very small, especially for lower degree eigenvalues below the so-called breakpoint $n_b$. The breakpoint is defined to be the value of the degree $l$ above which the radial functions $R_{ml}^{(2)}(-ic, i\xi)$ begin to increase in magnitude without bound while the corresponding radial functions $R_{ml}^{(1)}(-ic, i\xi)$ begin to decrease in magnitude. For small $m$ the breakpoint $n_b$ is approximately equal to $2(c_r + c_i)/\pi$, truncated to an integer. The Bouwkamp procedure often requires more accurate approximations or starting values for the lower order oblate eigenvalues than are available using traditional approximations and extrapolation. Oblfcn solves this problem by using a matrix method to obtain accurate values for these eigenvalues.

By successively choosing $n = 0, 1, 2, 3,...$, in (2), one can obtain an infinite set of simultaneous equations for the coefficients $d_n$. These equations can be written in matrix form as

$$\{B\}\{d\} = \lambda_{ml}\{d\}, \tag{4}$$

where $\{B\}$ is an infinite square tridiagonal matrix depending on $m$ and $c$, $\{d\}$ is a vector representation of the $d_n$ coefficients, and $\lambda_{ml}$ is the eigenvalue for $m$ and $l$. The desired oblate eigenvalues are then the set of eigenvalues of $\{B\}$.

Hanish and King [5] show that the matrix becomes symmetrical when $d_r$ is replaced with $d_r = \sqrt{\dfrac{(2r+2m+1)(r!)}{2(r+2m)!}}\, d_r$ , $r = n - 2, n$, and $n + 2$. This is equivalent to using associated Legendre functions with unit normalization in (1). It is much easier and faster to compute eigenvalues of a symmetric tridiagonal matrix. Furthermore, the $d_n$ coefficients with even subscript are only involved in spheroidal functions with even $l-m$, while those with odd subscript are involved when $l-m$ is odd. This allows the matrix $\{B\}$ to be decomposed into an even matrix $\{B^e\}$ using $n = 0, 2, 4,...$, whose eigenvalues are for $l-m$ even and an odd matrix $\{B^o\}$ using $n = 1, 3, 5,...$, whose eigenvalues are for $l-m$ odd. It is convenient to divide all of the matrix elements by $c^2$ and truncate both matrices to either order $4n_b/3$ or order 67, whichever is larger. Use of a standard tridiagonal matrix routine results in an odd and an even set of eigenvalues. Ordering each set of eigenvalues in increasing numerical value and interlacing the two sets results in accurate oblate eigenvalues $\lambda_{ml}(-ic)$ for $l - m = 0, 1, 2, ..., n_b$. The Bouwkamp procedure is attempted for all values of $l-m$ regardless of the value for $c$. This can provides eigenvalues that are slightly more accurate than the matrix results, as long as $l - m$ is not very small. The matrix values are used as starting values for $l-m$ up to $4n_b/3$ and estimates using extrapolation from previous eigenvalues are used for higher values of $l-m$. For lower values of $l-m$, when the Bouwkamp procedure fails to converge to an eigenvalue close to the matrix value, the matrix value is taken as the eigenvalue. When $c$ is large, neighboring low order eigenvalues beginning with $l = m$ and $l = m + 1$ are nearly identical.

For $c$ complex with a large $c_r$, the lower order eigenvalues are also paired including both the real and imaginary parts. When the magnitude of $c_i$ is small, the eigenvalues can still be ordered in increasing real part. When $c_i$ is greater than about 5, however, some of the intermediate eigenvalues are prolate like. The larger in magnitude that $c_i$ is, the greater the number of prolate-like eigenvalues. Usually these eigenvalues do not fit neatly in the eigenvalue sequence. The existence of prolate-like eigenvalues is apparently related to the fact that the recursion relation for the oblate angular function expansion coefficients can be obtained from the corresponding prolate recursion relation by replacing $c$ with $-ic$. Good estimates of these prolate-like eigenvalues are given by the standard asymptotic approximation [6, p. 243] for the lowest prolate eigenvalues with the oblate value for $c$ replaced by $-ic$, i.e., by $c_i - ic_r$. This approximation is:

$$\lambda_{ml}(-ic) \simeq -inc + m^2 + (n^2 + 5)/8 - in(n^2 + 11 - 32m^2)/(64c), \qquad (5)$$
$$where \quad n = 2(l-m) + 1.$$

The prolate-like eigenvalues are identified by their close numerical agreement with the asymptotic approximation. Coblfcn separately orders the eigenvalues for even and odd $l - m$ in increasing real part and combines the results. It then removes the prolate-like eigenvalues and places them in the sequence following either the eigenvalues with negative real part or additional paired eigenvalues with positive real part when they occur. This often provides a somewhat smooth transition between the prolate-like eigenvalues and the other eigenvalues. This method for ordering the eigenvalues is arbitrary, but it appears reasonable.

When $c$ is complex,, the matrix results are highly accurate only when $c_i$ is very small. They become increasingly inaccurate as $c_i$ increases, especially for very large $c_r$. Coblfcn uses an eigenvalue routine for complex tridiagonal matrices give by Cullum and Willoughby [7]. Conversion of the matrix routine used in oblfcn for real $c$ to complex arithmetic produced similar results so it appears that the matrices become less well-behaved as $c_i$ increases for large $c_r$. The least accurate eigenvalues are those near the breakpoint other than those that are prolate-like.

When $c_i$ is not very small and $c_r$ is moderate to large, the Bouwkamp procedure to refine the estimates can fail to converge to a very accurate eigenvalue. It can even fail to provide any improvement in accuracy beyond that obtained from the matrices. The problem is due to the nature of the numerator and denominator in the eigenvalue correction term [2]. The numerator is the difference between the value for $d_{l-m+2}/d_{l-m}$ obtained by forward recursion of (2) from $d_2/d_0$ or $d_3/d_1$, depending on whether $l - m$ is even or odd, and the value obtained by backward recursion from a sufficiently high value $n$ where $d_n/d_{n-2}$ is set equal to 0, The eigenvalue used in the recursion is the value after the previous iteration. After a few iterations, the two values for $d_{l-m+2}/d_{l-m}$ are essentially identical, agreeing to nearly all of the digits available in the precision used in the calculation. Continuing the Bouwkamp procedure cannot provide any further improvement. When $c$ is real, the correction term has a magnitude that is less than $10^{-ndec}$ times the eigenvalue, where ndec is the number of decimal digits used for real data. This results in an eigenvalue that is fully accurate, or nearly so. However, when $c$ is complex the correction term can have a relative magnitude that is nowhere near this small. Often the Bouwkamp procedure provides no improvement in the eigenvalue accuracy.

This is a bigger problem for double precision than for quadruple precision since the eigenvalues estimates from the matrices are often sufficiently accurate using quadruple precision. Coblfcn offers a solution to this problem when one is using double precision arithmetic and quadruple precision is available. One uses quadruple precision for the Bouwkamp procedure only. Now the Bouwkamp procedure can continue, if necessary, until the two values for $d_{l-m}/d_{l-m-2}$ agree to the number of decimal digits available in quadruple precision. This is sufficient to allow convergence of the eigenvalue to the number of digits in double precision, resulting in an eigenvalue that is fully accurate or nearly so. In addition to extending the range of parameters beyond those for which coblfcn provides useful results with double precision, it also improves the accuracy of the function values. And it does so with a modest increase in execution time. One can extend the useful parameter ranges for coblfcn even further by using quadruple precision for all calculations. But this extends the run time by a factor up to 50 or so.

An accurate eigenvalue allows one to now compute the $d$ coefficients using the recursion formula (2). Dividing each term in (2) by $d_n$ results in an expression relating the ratio $N_{n+2} = d_{n+2}/d_n$ to the ratio $N_n = d_n/d_{n-2}$. Traditionally this expression is used in the forward direction to calculate ratios up to $N_{l-m}$ starting with the first ratio $N_2 = (\beta_0 - \lambda_{ml})/\alpha_0$ for $n$ even or $N_3 = (\beta_1 - \lambda_{ml})/\alpha_1$ for $n$ odd. Ratios for $n$ above $l-m$ are calculated backward from a suitably high value of $n$ where the ratio $N_{n+2}$ is set equal to zero. As the ratios are calculated backward, they become progressively more accurate until they are essentially fully accurate. Only the backward recursion is used when $l-m = 0$ or 1. See Flammer [2] for a discussion of this. This procedure always works well for the prolate case. Here the forward recursion provides accurate values up the ratio $N_{l-m}$ and the backward recursion provide accurate values down to the ratio $N_{l-m+2}$. This is not always true in the oblate case when $c$ is greater than about 50.

Reference [1] describes a procedure that provide accurate values for the ratios of the oblate $d$ coefficients. One first uses the backward recursion down to the either $N_2$ for $n$ even or $N_3$ for $n$ odd. Then one uses the forward recursion from either $N_2$ or $N_3$ until the forward and backward ratio values match to ndec decimal digits, where ndec is the number of decimal digits available in real numbers. If there is no match to ndec digits, the forward recursion is continued until the match starts decreasing significantly. The best match is then selected. It is rarely less than ndec - 2 digits. This procedure works equally well for complex $c$ where the match is based on the magnitude of the two ratio values. However, when the eigenvalue is not fully accurate or nearly so, the match is typically to the number of accurate digits in the eigenvalue. This procedure is used in coblfcn when $c_r$ is greater than 50.

The resulting coefficients can be normalized by requiring that $S_{ml}^{(1)}(-ic,\eta)$ has the same normalization factor as $P_l^m(\eta)$ [5], resulting in the following relation:

$$\sum_{n=0,1}^{\infty} {}' \frac{2(n+2m)!}{[2(n+m)+1]n!}[d_n(-ic\,|\,ml)]^2 = \frac{2(l+m)!}{(2l+1)(l-m)!}. \tag{6}$$

Use of this Meixner and Schäfke [6] normalization scheme has the practical advantage of eliminating the need to compute the normalization factor which is often found in problems involving expansions in spheroidal angular functions. A better choice, however, is to set the rhs of (6) equal to unity. This results in the angular functions having unit norm. It has the advantage of limiting the magnitude of the angular functions to moderate values. For other normalization schemes such as (6) the angular functions become increasing large in magnitude as *m* increases and eventually overflow. Both oblfcn and coblfcn offer either unit normalization or the same normalization as the corresponding Legendre functions.

When *c* is real, the normalization sum in (6) is numerically robust with no subtraction errors occurring in its computation. This contrasts with the corresponding Flammer normalization sum involved in requiring the angular functions to match the corresponding associated Legendre functions at $\eta = 0$ [see e.g., Flammer [2, p. 21]. Here subtraction errors can occur in calculating the angular function at $\eta = 0$, especially when *c* is large and $l - m$ is less than the breakpoint $n_b$. Subtraction error is defined to be the number of accurate decimal digits that are lost in calculating the sum of the series. This loss of accuracy occurs when the sum of all of the positive terms in the series is nearly equal to the sum of all of the negative terms. The subtraction error is then equal to the number of leading decimal digits that are the same in the positive and negative sums. When *c* is complex, subtraction error can also occur in the sum in (6). The errors are maximum for *m* = 0. They are small for low values of *l*, increase to a maximum at a value of *l* somewhat below the breakpoint, and then decrease as *l* increases further. The maximum error is no more than about 1 digit when $c_i$ is less than 10 but increases rapidly as $c_i$ increases. Figure 1 shows the maximum subtraction error in decimal digits at *m* = 0 plotted versus $c_i$ for selected values of $c_r$.

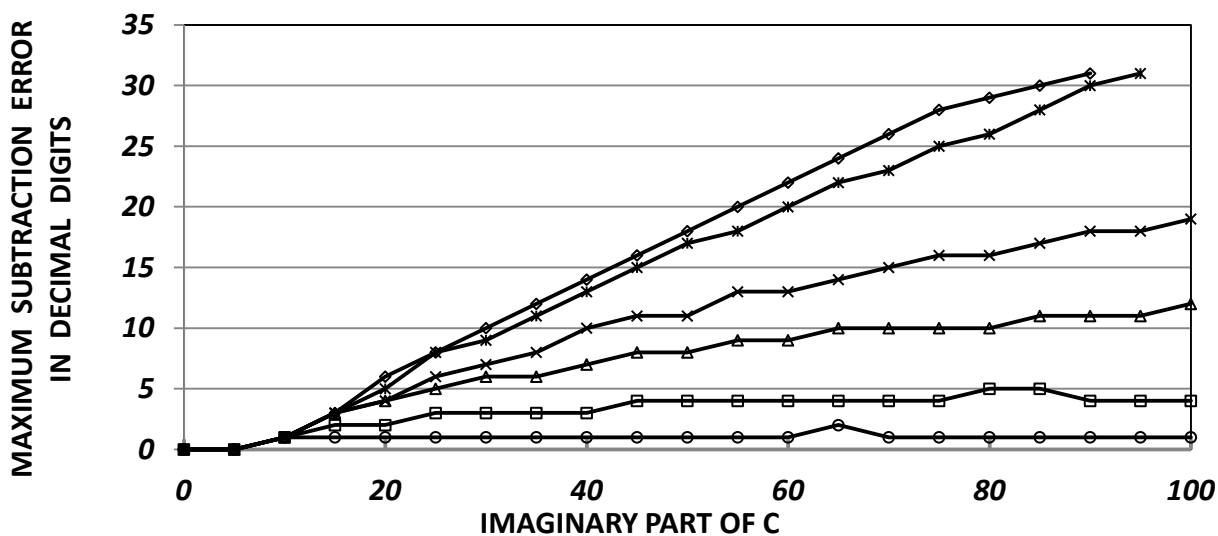

Fig 1: Maximum subtraction error in decimal digits involved in calculating the Meixner and Schäfke normalization plotted versus the imaginary part of c for selected values of the real part of c: o (10), □ (20), ∆ (50), × (100), * (500), ◊ (1000).

Calculation of the angular functions using (1) can also suffer subtraction errors at large values of *c* for values of $\eta$ other than 0. For a given value of *m* and for $l - m$ less than $n_b$, the error decreases to zero as $\eta$ increases from 0 to 1. For given values of *m* and $\eta$, the error decreases to zero as $l - m$ increases. Reference [1] shows a graph illustrating this behavior when *c* is real. Similar behavior occurs when *c* is complex. When subtraction error is encountered in the computation of (1), the resulting angular functions and their first derivatives are reduced in magnitude by an amount corresponding to the subtraction error. Their magnitude in this case is corresponding smaller than angular functions for higher values of *l* and/or $\eta$ not near zero. The loss in accuracy due to these subtraction errors will not likely affect numerical results for physical problems using these functions.

A third source of inaccuracy in the angular functions arises from the potential loss of accuracy in the eigenvalues at values of *l* near and somewhat below the breakpoint when $c_i$ is not very small and $c_r$ is moderate to large. As discussed above, eigenvalues that are fully accurate or nearly so are obtained for double precision by using quadruple precision for the Bouwkamp procedure to ensure full convergence. Otherwise, the eigenvalue accuracy is estimated using the degree of convergence of the Bouwkamp procedure as well as the degree of pairing of neighboring eigenvalues when applicable. One would normally expect that a reduction in the accuracy of the eigenvalue would result in a corresponding decrease in the accuracy of the Meixner and Schäfke normalization beyond that resulting from subtraction errors in their calculation. It turns out to be more complicated than this. Comparison of double precision results, double precision results using quadruple precision for the Bouwkamp procedure, and quadruple precision results show that a good estimate of the accuracy of the Meixner and Schäfke normalization is given by the smaller of naccre - 1, itestm - 1, and ndec - jsubms - 1. Here naccre is the estimated accuracy of the eigenvalue, itestm is the number of digits of agreement between the forward and backward recursions to determine the *d* coefficients, ndec is the precision in the arithmetic used, and jsubms is the subtraction error in the Meixner and Schäfke normalization. The same effect is seen regarding the loss in accuracy of the angular functions due to subtraction errors in evaluating (1).

The resultant accuracy of the angular functions $S_{ml}^{(1)}(-ic,\eta)$ is conservatively estimated using the subtraction error involved in evaluating the series in (1), the subtraction error involved in calculating the Meixner and Schäfke normalization, the estimated accuracy of the eigenvalue and the number of decimal digits that match in the forward and backward recursions to compute the *d* coefficients.

## 3      Expansion of the product of the radial and angular functions

The expansion of the product of $R_{ml}^{(j)}(-ic,i\xi)$ and $S_{ml}^{(1)}(-ic,\eta)$ in terms of the corresponding spherical functions is given by:

$$R_{ml}^{(j)}(-ic,i\xi) S_{ml}^{(1)}(-ic,\eta) = \sum_{n=0,1}^{\infty} {}' i^{n+m-l} d_n(-ic|ml) \psi_{n+m}^{(j)}(kr) P_{n+m}^{m}(\cos\theta), \qquad (7)$$

where *j* = 1 or 2. $\psi_{m+n}^{(1)}(kr)$ is the spherical Bessel function $j_{m+n}(kr)$ and $\psi_{m+n}^{(2)}(kr)$ is the spherical Neumann function $y_{n+m}(kr)$. Here *c* can be complex. This is a special case of the more general expansion given by Meixner and Schäfke [6, p. 307]. Using the relationship between the spherical coordinates *r* and $\theta$ and spheroidal coordinates (about the same origin and with $\eta$ = 1 coincident with $\theta$ = 0) we obtain $kr = c(\xi^2 - \eta^2 + 1)^{1/2}$ and $\cos\theta = \eta\xi/(\xi^2 - \eta^2 + 1)^{1/2}$. Substituting for $S_{ml}^{(1)}(-ic,\eta)$ from (1) and solving for $R_{ml}^{(j)}(-ic,i\xi)$ produces

$$R_{ml}^{(j)}(-ic,i\xi) = \frac{\sum_{n=0,1}^{\infty} {}' i^{n+m-l} d_n(-ic|ml) \psi_{n+m}^{(j)}[c(\xi^2-\eta^2+1)^{1/2}] P_{n+m}^m[\eta\xi/(\xi^2-\eta^2+1)^{1/2}]}{\sum_{n=0,1}^{\infty} {}' d_n(-ic|ml) P_{n+m}^m(\eta)}. \qquad (8)$$

The significance of this general expression is that it allows us to choose the value for $\eta$ that provides the maximum accuracy for calculated values of $R_{ml}^{(j)}(-ic,i\xi)$. For many parameter ranges it is desirable to allow $\eta$ to vary as the value of the index $l$ increases from $m$ to higher values. Reference [1] describes the application of (8) to the calculation of the radial functions of both the first and second kinds for real values of $c$. This paper will address the use of (8) when $c$ is complex.

## 4  Traditional Bessel function expressions

Consider the case when $\eta = 1$. The argument of $P_{n+m}^m$ in both the numerator and the denominator approaches unity as $\eta$ approaches unity. Although $P_{n+m}^m$ approaches zero in this case for $m \neq 0$, the limit of the rhs of (8) exists and we obtain:

$$R_{ml}^{(j)}(-ic,i\xi) = \left(\frac{\xi^2+1}{\xi^2}\right)^{m/2} \frac{\sum_{n=0,1}^{\infty} {}' i^{n+m-l} d_n(-ic|ml) \psi_{m+n}^{(j)}(c\xi) \frac{(n+2m)!}{n!}}{\sum_{n=0,1}^{\infty} {}' d_n(-ic|ml) \frac{(n+2m)!}{n!}}. \qquad (9)$$

Flammer [2, p. 32] derives (9) using integral representations of the spheroidal wave functions. The corresponding expression for the first derivative of $R_{ml}^{(j)}$ with respect to $\xi$ is obtained by taking the first derivative of the rhs of (9). Equation (9) is the expression commonly used to calculate numerical values for both $R_{ml}^{(1)}$ and $R_{ml}^{(2)}$. The advantage of these expressions for the oblate case when $c$ is real is that the denominator sum is robust with no subtraction errors. This sum is the one involved in the Morse and Feshbach scheme for normalizing the angular functions so that they are equal at $\eta = 1$ to the corresponding associated Legendre function at $\cos\theta = 1$. When $c$ is complex, significant subtraction errors can occur in evaluating the denominator sum. For a given value of $c_r$, $c_i$, and $m$, the error is a maximum at the degree $l$ where the first prolate-like eigenvalue occurs. The subtraction error decreases rapidly as $l$ increases or decreases from this point. The maximum subtraction error is essentially independent of $c_r$. It is not surprising that the maximum subtraction error is equal to the error in computing $S_{mm}^{(1)}(c,1)$ for prolate functions with $c$ real and equal to $c_i$ [8, Fig.1] since the lowest prolate-like eigenvalue corresponds to $l = m$. Figure 2 shows the maximum subtraction error in decimal digits plotted versus $c_i$ for selected values of $m$.

## 5  Alternative Bessel function expressions

Reference [2] also provides expressions for the prolate radial functions of the first kind obtained by choosing $\eta = 0$ in the prolate version of (8). Converting these expressions to oblate form and extending them to include radial functions of the second kind results in the following alternative expressions:

$$R_{ml}^{(j)}(-ic,i\xi) = \frac{\sum_{n=0,1}^{\infty}{}' i^{n+m-l} d_n(-ic|ml) \psi_{m+n}^{(j)}(c\sqrt{\xi^2+1}) P_{n+m}^m(0)}{\sum_{n=0,1}^{\infty}{}' d_n(-ic|ml) P_{n+m}^m(0)}, \quad l-m, \text{ even} \quad (10)$$

$$R_{ml}^{(j)}(-ic,i\xi) = \left(\frac{\xi}{\sqrt{\xi^2+1}}\right) \frac{\sum_{n=0,1}^{\infty}{}' i^{n+m-l} d_n(-ic|ml) \psi_{m+n}^{(j)}(c\sqrt{\xi^2+1}) \left[\frac{dP_{n+m}^m(\eta)}{d\eta}\right]_{\eta=0}}{\sum_{n=0,1}^{\infty}{}' d_n(-ic|ml) \left[\frac{dP_{n+m}^m(\eta)}{d\eta}\right]_{\eta=0}}, \quad l-m \text{ odd.} \quad (11)$$

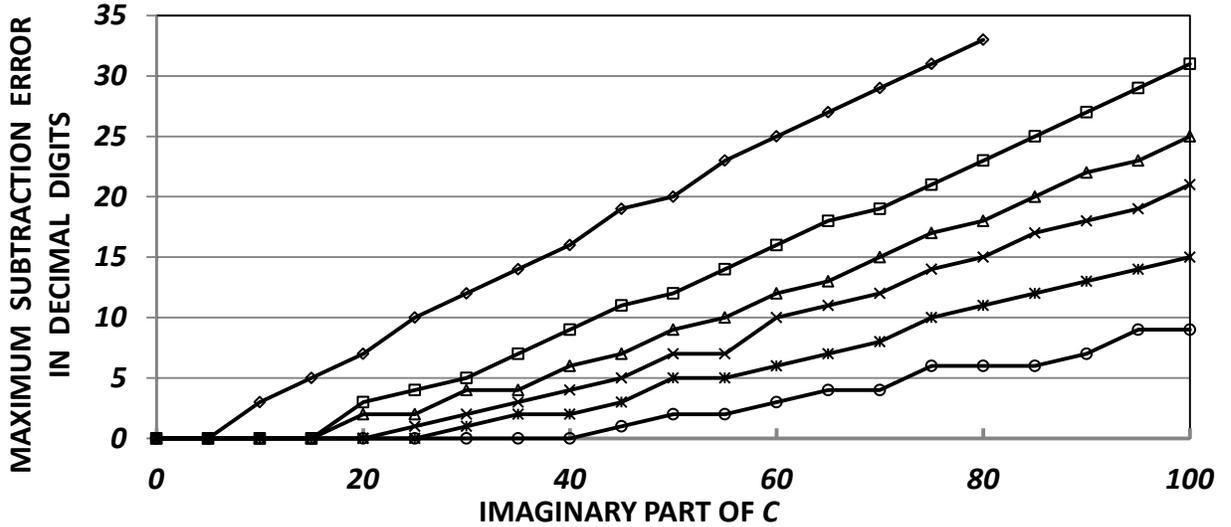

Fig 2: Maximum subtraction error in decimal digits involved in calculating the Morse and Feshbach normalization plotted versus the imaginary part of $c$ for selected orders $m$: ◊ (0), □ (20), Δ (40), × (60), * (100), o (200).

Reference [8] shows that the corresponding expressions for the prolate case are numerically robust when calculating the prolate radial function of the first kind. Neither the numerator nor the denominator sums suffer subtraction errors so that accurate values are obtained over all parameter ranges. For the oblate case this is not true, even for real $c$. The denominators in (10) and (11) are the Flammer normalization sums for $l - m$ even and $l - m$ odd, respectively. The $\eta = 0$ expressions will prove useful in the calculation of the radial functions of the second kind.

# 6    Calculation of the oblate radial functions of the first kind

The traditional expression (9) provides highly accurate values for the oblate radial functions of the first kind when $c$ is small or $l$ is somewhat larger than the break point $n_b$. When $c$ is complex, both $c_r$ and $c_i$ must be small. There are subtraction errors in calculating the numerator term in (9) that increase with increasing $c_r$. These errors are maximum at $l - m = 0$ and then decrease monotonically to zero with increasing $l - m$. Figure 4 in [1] shows examples of the subtraction error arising when calculating $R_{ml}^{(1)}$ and its first derivative using the traditional Bessel function expression (9) for real values of $c$. The subtraction errors are zero or nearly so for $m = 0$, increase with increasing $m$ until $m$ is very large and then decrease with increasing $m$. They also tend to be largest when $\xi$ is near unity. Moreover, the subtraction errors plotted in Fig. 4 in [1] occur in calculating the numerator terms, since the denominator does not incur subtraction errors when $c$ is real. As seen above in Fig. 2, there can be significant subtraction errors arising in calculation of the denominator term when $c$ is complex. The corresponding subtraction errors arising from the numerator terms for $c$ complex are smaller than those for real $c$ and decrease as $c_i$ increases in value. Examples of this are given in Fig. 3. Here the subtraction errors are plotted as a function of $l - m$ both for real $c$ and for the case where $c_i = 20$. The subtraction errors continue to decrease relative to those for real $c$ as $c_i$ increases above 20. The subtraction errors shown in Figure 3 are for $R_{ml}^{(1)}$. The corresponding errors for its first derivative are nearly identical. Note that the subtraction error that limits the accuracy of (9) is the larger of the error for the numerator and the denominator.

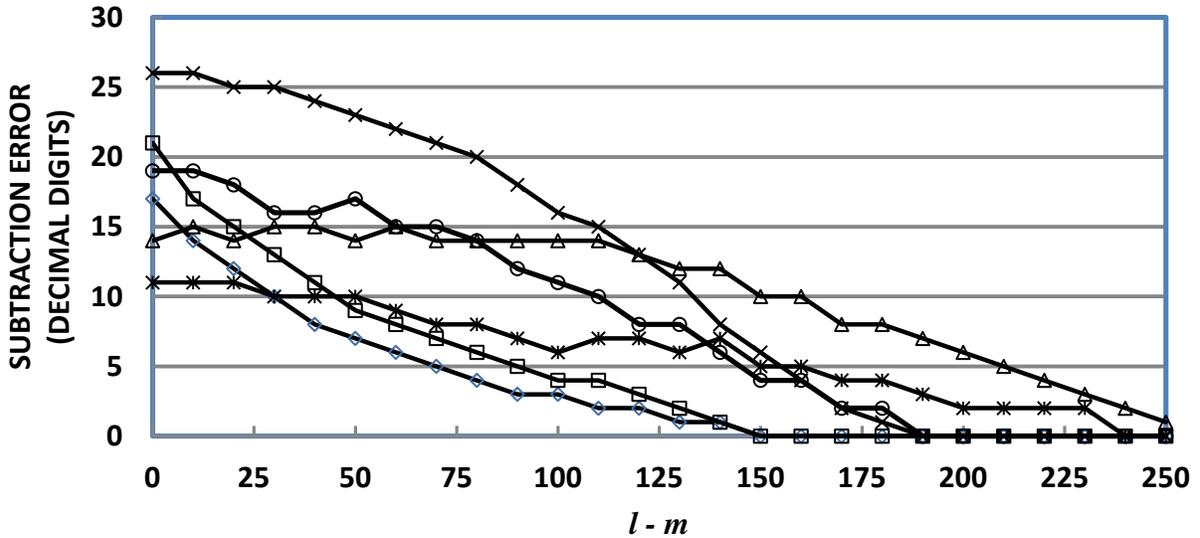

Fig. 3: Subtraction error when calculating the numerator of $R_{ml}^{(1)}(-ic, i\xi)$ using the traditional Bessel function expression for selected parameters ($\xi$, $c_r$, $c_i$, $m$): □ (0.5, 300, 0, 100); ◊ (0.5, 300, 20, 100); Δ (1.0, 300, 0, 100); * (1.0, 300, 20, 100); × (1.0, 300, 0, 200); o (1.0, 300, 20, 200).

These subtraction errors can be avoided by using the variable $\eta$ method to calculate $R_{ml}^{(1)}$ and its first derivative. If the subtraction errors for $l = m$ are not near zero using the traditional expressions, then the value for $\eta$ is reduced in steps from unity until the errors are as small as possible. Usually both subtraction errors will now be near zero. However, if either $R_{ml}^{(1)}$ or its first derivative is near a root, its value will be somewhat smaller in magnitude than expected and its calculation will involve an unavoidable subtraction error whose size depends on how close the root is.

Using equal steps in $\theta = \arccos(\eta)$ works well, with a step size $\Delta\theta$ of about 0.05 radian. Here $\theta$ is incremented from zero, the value for the traditional expressions. Only a few steps are usually required. The value of $\theta$ that worked for $l = m$ is then used for progressively higher values of $l$ until the subtraction error increases. Then the value of $\theta$ is incrementally decreased until the subtraction error is near zero again. Usually only one step is required. The process is continued either until the maximum desired value of $l$ is reached or until $\theta$ has reached zero. Once $\theta$ reaches zero, the traditional expressions work well for higher values of $l$.

When $c$ is complex, the process becomes more complicated. The search for the best value of $\eta$ often takes more than one step when both $c_r$ and $c_i$ are large and when near the breakpoint. It is sometimes difficult in this case to achieve subtraction errors that are near zero. Coblfcn strives to achieve errors no greater than 2 digits when using double precision and no greater than 4 digits when running quadruple precision for those cases where double precision is inadequate. Care has to be taken when reaching the degree $l$ where the first prolate-like eigenvalue occurs. Also, sometimes one goes well beyond the breakpoint before the traditional expression works well with small subtraction errors.

The accuracy of $R_{ml}^{(1)}$ and its first derivative is estimated using the subtraction errors that occurs in their series calculation including the denominator series, the estimated accuracy of the eigenvalue, and the number of decimal digits of agreement between the forward and backward recursions to compute the $d$ coefficients.

## 7  Use of low-order pairing of eigenvalues at higher values of $c$

Calculation of the radial functions of the second kind can be difficult even when $c$ is real. This is especially true at low orders when $c_r$ is large and $\xi$ is somewhat less than unity. Traditional expansions often fail here. Fortunately one can take advantage of the low-order pairing of eigenvalues. When $c_r$ is not small, the eigenvalue $A_{mm}$ is nearly equal to $A_{m,m+1}$; $A_{m,m+2}$ is nearly equal to $A_{m,m+3}$; .... The agreement decreases with increasing $l$ and disappears as $l - m$ approaches the breakpoint $n_b$. The differential equations for the corresponding low-order radial functions are then nearly identical as are their solutions. The following approximations are found to apply here:

$$R_{ml}^{(2)}(-ic, i\xi) \simeq R_{m,l+1}^{(1)}(-ic, i\xi), l - m \text{ even},$$
$$\simeq -R_{m,l-1}^{(1)}(-ic, i\xi), l - m \text{ odd}. \tag{12}$$

This is equivalent to setting $R_{m,m+1}^{(3)} = R_{m,m+1}^{(1)} + iR_{m,m+1}^{(2)} \simeq -i\left[ R_{m,m}^{(1)} + iR_{m,m}^{(2)} \right]$ and so forth for higher values of $l$ with paired eigenvalues. First derivatives of the radial functions also satisfy the same relationships. Note that there is a typo in this equation in [1] where $l$ even and $l$ odd appears rather than the correct $l - m$ even and $l - m$ odd. Flammer [2, p.68] obtained these approximations for the special case of $\xi = 0$ based on an asymptotic representation of the radial function $R_{m,l}^{(3)}$ in a series of Laguerre functions. Observation of numerical results shows that (12) is equally valid for all values of $\xi$. It also shows that the number of decimal digits of agreement predicted by (12) tends to be the number of digits of agreement between the corresponding eigenvalues, with an exception when $\xi$ is very small. Here, when $l - m$ is even, the number of digits of agreement between $R_{ml}^{(2)}$ and $R_{m,l+1}^{(1)}$ tends to be smaller than this by the logarithm to the base 10 of $\xi$, truncated to an integer. The same is true for the first derivatives of $R_{ml}^{(2)}$ and $-R_{m,l-1}^{(1)}$ when $l - m$ is odd. When well below the breakpoint, both $R_{ml}^{(2)}$ for $l - m$ even and its first derivative for $l - m$ odd have magnitudes that become increasingly small with decreasing $\xi$. The lower accuracy of these function values obtained using (12) may not be a problem since it is expected that their contributions to the solution of problems involving these functions would also become increasingly small.

Reference [1] includes a figure that shows an example of the near equality between the number of decimal digits of agreement between neighboring eigenvalues and the number of digits of agreement between the radial functions in accordance with (11), This example is for real $c$, but similar results are obtained for $c$ complex. Coblfcn uses (12) to obtain values for $R_{ml}^{(2)}$ and its first derivative when applicable. Care is taken to account for the lower accuracy expected when $\xi$ is very small so that coblfcn can try other methods and possibly obtain more accurate values than provided using (12).

When using (12), coblfcn does not estimate the accuracy based on the agreement between neighboring eigenvalues but rather uses the Wronskian relationship:

$$R_{ml}^{(1)} \frac{dR_{ml}^{(2)}}{d\xi} - R_{ml}^{(2)} \frac{dR_{ml}^{(1)}}{d\xi} = \frac{1}{c(\xi^2 + 1)}. \tag{13}$$

An integer estimate of accuracy is given by the number of leading digits of agreement between the calculated Wronskian on the lhs of (13) and the theoretical Wronskian given by the rhs side. The traditional Bessel function expressions together with the variable $\eta$ method almost always provide accurate values for $R_{ml}^{(1)}$ and its first derivative. Then the number of decimal digits of agreement between the theoretical and computed Wronskian is a measure of the accuracy of $R_{ml}^{(2)}$ and its first derivative. However, when $\xi$ is very small, one of the two products on the lhs of (13) tends to be much smaller than the other. Use of the Wronskian can then overestimate the accuracy of $R_{ml}^{(2)}$ for $l - m$ even and its first derivative for $l - m$ odd. Usually $R_{ml}^{(2)}$ and its first derivative are calculated with nearly the same accuracy. However, this is not true for very small values of $\xi$ when using either (12) or the integral method that will be discussed below in Sec. 11.

## 8   Use of near equality of $iR_{ml}^{(1)}(c\xi)$ and $R_{ml}^{(2)}(c\xi)$ for large values of $c_i\xi$

When the product of the radial coordinate $\xi$ and $c_i$ increases to large values, $R_{ml}^{(1)}$, $R_{ml}^{(2)}$, and their first derivatives can become large in magnitude for lower values of $l - m$. This same behavior is seen in the Bessel functions that appear in the traditional expansions for the radial functions. Here

$$j_0(c\xi) = \sin(c\xi)/c\xi = (e^{ic\xi} - e^{-ic\xi})/2ic\xi$$
$$= -i(e^{ic_r\xi - c_i\xi} - e^{-ic_r\xi + c_i\xi})/2c\xi \simeq -ie^{-ic_r\xi + c_i\xi}/2c\xi, \text{ as } c_i\xi \to \infty, \quad (14)$$

and

$$y_0(c\xi) = -\cos(c\xi)/c\xi = -(e^{ic\xi} + e^{-ic\xi})/2c\xi$$
$$= -(e^{ic_r\xi - c_i\xi} + e^{-ic_r\xi + c_i\xi})/2c\xi \simeq -e^{-ic_r\xi + c_i\xi}/2c\xi, \text{ as } c_i\xi \to \infty. \quad (15)$$

Thus as $c_i\xi$ becomes large, both $j_0(c\xi)$ and $y_0(c\xi)$ have magnitudes nearly equal to $e^{c_i\xi}/2c_m\xi$, where $c_m$ is the magnitude of $c$. Also $y_0(c\xi)$ becomes very nearly equal to $ij_0(c\xi)$. As the order increases, this relationship tends to persist but the near equality decreases as both $j_n(c\xi)$ and $y_n(c\xi)$ decrease in magnitude to a value approximately equal to $1/c_m\xi$,. Note that it has been assumed that $c_i$ is positive. If $c_i$ is negative, both have magnitudes nearly equal to $e^{-c_i\xi}/2c_m\xi$

When the corresponding radial functions become large, they satisfy the approximation $R_{ml}^{(2)}(c,\xi) \simeq iR_{ml}^{(1)}(c,\xi)$ and similarly for their first derivatives. The agreement between $R_{ml}^{(1)}$ and $iR_{ml}^{(2)}$ in decimal digits is given by the quantity acc = int{ $\log 10(abs[R_{ml}^{(1)}(dR_{ml}^{(1)}/d\xi)c(\xi^2 + 1)])$ }. Here abs refers to absolute value and int denotes truncation to an integer. This agreement is required in order for the Wronskian relationship (13) to be satisfied. The integer acc is determined after $R_{ml}^{(1)}$ and its first derivative have been accurately calculated. The value for $R_{ml}^{(2)}$ accurate to acc digits or to the estimated accuracy of $R_{ml}^{(1)}$, whichever is smaller, is given by $iR_{ml}^{(1)}$. Similarly for the first derivative of $R_{ml}^{(2)}$.

For $c_i\xi$ not small, the radial function of the first kind with the largest magnitude corresponds to the lowest order prolate-like eigenvalue. If $c_r$ is very small, this will be located at $l = m$. For larger values of $c_r$, the lowest prolate-like eigenvalue is located immediatly after all of the eigenvalues with negative real parts and possible paired eigenvalues with positive real parts. The magnitude of the largest radial function of the first kind tends to be somewhat less than $e^{c_i\xi}/2c_m\xi$.

One must take care in using the Wronskian to estimate the accuracy of calculated radial functions when $R_{ml}^{(1)}$ has a magnitude somewhat larger than $1/c_m\xi$. A subtraction error of acc digits will occur when computing the lhs of (13), assuming the function values are accurate to at least this many digits. The comparison with the theoretical value given by the rhs of (13) will thus be reduced by acc digits. A good accuracy estimate is given by adding acc digits to the Wronskian comparison result. However, coblfcn does not automatically use this estimate. When acc is at least one, coblfcn also estimates the accuracy using subtraction errors involved in the calculation of both $R_{ml}^{(2)}$ and its first derivative as well as the estimated accuracy of $R_{ml}^{(1)}$ and its first derivative. The smaller of these two estimates is used. When the Legendre function expression is used, acc is only added to the Wronskian comparison result if acc exceeds one.

## 9    Traditional Legendre function expression

A variety of methods are required to obtain accurate values for $R_{ml}^{(2)}$ over a wide range of parameters. Considered next is use of the traditional Legendre function expression {see, e.g.,

Flammer [2, p. 33] or Zhang and Jin [9, p. 567]}. Reference [10] shows that this expression can be obtained using the product expansion. Here the angular function of the first kind is replaced with that of the second kind and the associated Legendre function of the first kind $P_{n+m}^m$ is replaced with the associated Legendre function of the second kind $Q_{n+m}^m$. The resulting expression is:

$$R_{ml}^{(2)}(-ic, i\xi) = \frac{1}{\kappa_{ml}^{(2)}(-ic)} \sum_{n=-\infty}^{\infty}{}' d_n(-ic\,|\,ml) Q_{n+m}^m(i\xi). \qquad (16)$$

Of course, equation (16) could have been obtained more directly from the fact that $R_{ml}^{(2)}$ and $S_{ml}^{(2)}$ are proportional to each other.

Although $Q_{n+m}^m(i\xi)$ becomes infinitely large when $n$ is less than $-2m$, its product with $d_n$ is finite and proportional to $P_{-n-m-1}^m(i\xi)$. The rhs of (16) then divides into two series, one over $n$ from $-2m$ (or $-2m + 1$ if $l - m$ is odd) to $\infty$ involving $Q_{n+m}^m(i\xi)$ and one over $n$ from $2m + 2$ (or $2m + 1$ if $l - m$ is odd) to $\infty$ involving $P_{n-m-1}^m(i\xi)$. The result is the traditional associated Legendre function expression used to evaluate $R_{ml}^{(2)}(-ic, i\xi)$ when $\xi$ is small.

The $d_n$ coefficients with negative subscripts required in the $Q_{n+m}^m(i\xi)$ series and the special $d_{\rho/n}$ coefficients required for the $P_{n-m-1}^m(i\xi)$ series can be computed from either $d_0$ or $d_1$, depending on whether $l - m$ is even or odd. The expression for the joining factor $\kappa_{ml}^{(2)}(-ic)$ given by Flammer [2] contains the same series that is given in the denominator of the traditional Bessel function expression (9). For $c$ real this series does not suffer subtraction errors but does so for $c$ complex as shown in Fig. 2. It also explicitly contains the Flammer normalization sum due to the presence of $d_{-2m}$ or $d_{-2m+1}$ in the expression. This sum is basically the denominator shown above in (10) and (11). Figure 2 in Reference [1] illustrates that for real $c$, the Flammer normalization sum suffers a loss in accuracy for low values of $l - m$ due to subtraction errors that increase without bound as $c$ increases. For a given value of $c$, the subtraction errors are a maximum at $l = m$ and decrease to zero as $l$ increases. Note that the labels for Fig. 2 in [1] are incorrect. The 6 plots shown there represent m = 0, 10, 20, 30, 40, and 50. When $c$ is complex, the maximum error at $l = m$ is nearly the same as shown in Fig. 2 from [1], regardless of the value for $c_i$. However, the subtraction error decreases more slowly with increasing $l$ as $c_i$ increases. When the first prolate-like eigenvalue is reached, the subtraction error drops to 0 and then slowly increases until is it back to the value it had immediately before the first prolate-like eigenvalue. After this, it slowly decreases to zero as $l$ continues to increase.

Another source of inaccuracy in the joining factor can arise from the factor $d_{-2m}$ for $l - m$ even and $d_{-2m+1}$ for $l - m$ odd. Values for the ratios of successive $d_n$ coefficients with negative index $n$ are obtained using backward recursion on (2). One starts with the closed form expression for $d_{-2m}/d_{-2m+2}$ when $l - m$ is even and continues to $d_{-2}/d_0$. For $l - m$ odd, one starts with $d_{-2m+1}/d_{-2m+3}$ and ends with $d_{-1}/d_1$. The ratio $d_{-2m}/d_0$ or $d_{-2m+1}/d_1$ is calculated by multiplication of successive coefficient ratios. When $c$ is real and large or $c_r$ is large and $l$ is below the breakpoint, significant subtraction errors occur during the recursion process. This reduces the accuracy of the joining factor as well as the sum involving $Q_{n+m}^m(i\xi)$ with negative values appearing in (16). The subtraction error in calculating $d_{-2m}$ and $d_{-2m+1}$ is determined in coblfcn by examination of the recursion process. Figure 4 shows this subtraction error plotted as a function of $l - m$ for selected parameters. The dips in the curves for $c_i = 20$ and 40 occur where the prolate-like eigenvalues are located.

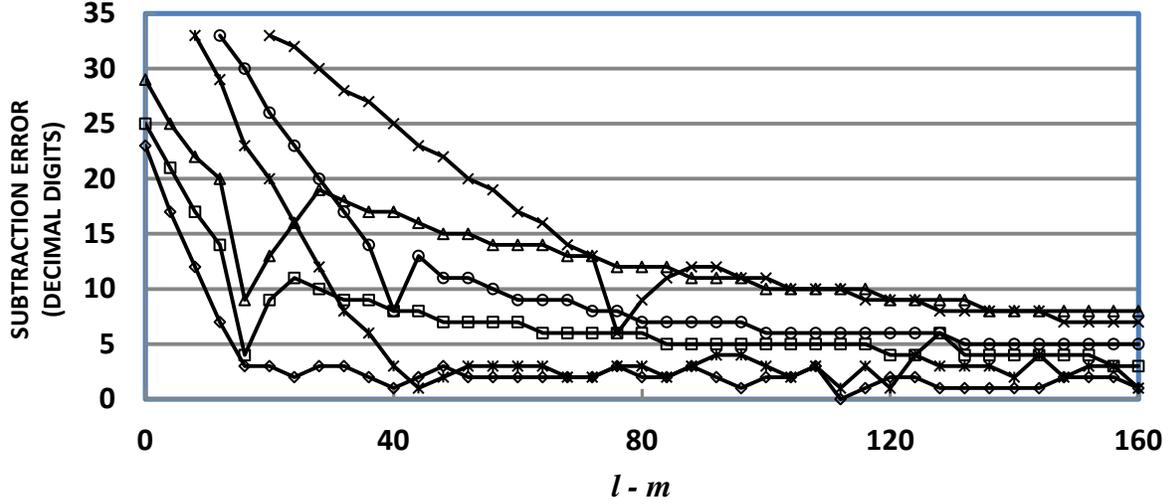

Fig. 4: Subtraction error when calculating $d_{-2m}$ and $d_{-2m+1}$ for selected parameters $(c_r, c_i, m)$: ◊ (100, 0, 50); □ (100,20, 50); Δ (100, 40, 50); * (200, 0, 100); × (200,40, 50); o (200,20, 100).

It is convenient to compute all of the $d$ coefficients as ratios of successive coefficients, Consider the $d_{\rho/n}$ coefficients. One calculates them by forward recursion using (2) starting with a sufficiently high value of $n$ where the ratio is set equal to 0. The process continues until the final step produces the ratio $d_{\rho/2m+2}/d_{-2m}$ or $d_{\rho/2m+1}/d_{-2m+1}$, depending on whether $l - m$ is even or odd. When $c$ is real or $c_i$ is very small, a subtraction error can occur just during this step, resulting in reduced accuracy for the lead coefficient in both the $P^m_{n-m-1}(i\xi)$ sum and the corresponding sum for the first derivative of $R^{(2)}_{ml}$. The error can be large enough to preclude radial functions of the desired accuracy. After the publication of [1], a procedure was added to oblfcn that can sometimes improve the accuracy of this coefficient. Here the Wronskian relationship is solved for the coefficient in terms of the calculated values for each of the sums in the Legendre expression for both $R^{(2)}_{ml}$ and its first derivative. The resulting new value for the coefficient is then used whenever the accuracy is expected to improve. An estimate of the improved accuracy is calculated using subtraction errors in the various sums together with the estimated accuracy of $R^{(1)}_{ml}$ and its first derivative. Comparison of double precision and quadruple precision results show this estimate to be reliable. The procedure is used when $\xi \leq 0.01$, primarily near the break point. Here it can sometimes provide more accurate values for $R^{(2)}_{ml}$ when $l - m$ is even and for its first derivative when $l - m$ is odd than those given by eigenvalue pairing or the integral method. Note that when $c_i$ is not very small, subtraction errors can also occur in the calculation of many of the $d_{\rho/n}$ coefficients ratios, rather than just the last one. However, there is often a large subtraction error in the last step, which allows the procedure.

Legendre function sums in (16) can also suffer subtraction errors that are as large as occur in the joining factor, especially for low values of $m$. Figure 7 in [1] shows examples of the subtraction error in the Legendre function sums for selected parameter sets when $c$ is real. When $c$ is complex, the behavior is somewhat similar except the subtraction error tends not to decrease as rapidly with increasing $l - m$ as shown in [1]. Also it can be larger for higher values of $m$. However, unless $c_i$ is small or $m = 0$, the errors in the joining factor arising from the $d$ coefficients with negative index are often much larger than the subtraction errors in the Legendre function sums. For $c$ complex, the traditional Legendre function expression is used when $\xi$ is less than or equal to 0.99, just as for $c$ real. Other methods will be used to compute $R^{(2)}_{ml}$ and its first derivative when $\xi$ is larger or when the results from the traditional Legendre function expression are not sufficiently accurate.

Examples of the accuracy obtained using (16) are given in Fig. 5. Calculations were carried out in quadruple precision with 33 decimal digits of precision. The accuracy is plotted versus $l - m$ for $\xi = 0.1$ with $m = 100$ and for $\xi = 0.99$ with $m = 50$. Results are shown for $c_i = 0$, 20, and 40. The accuracy for small values of $m$ varies much less with $c_i$ than shown in Fig. 5.

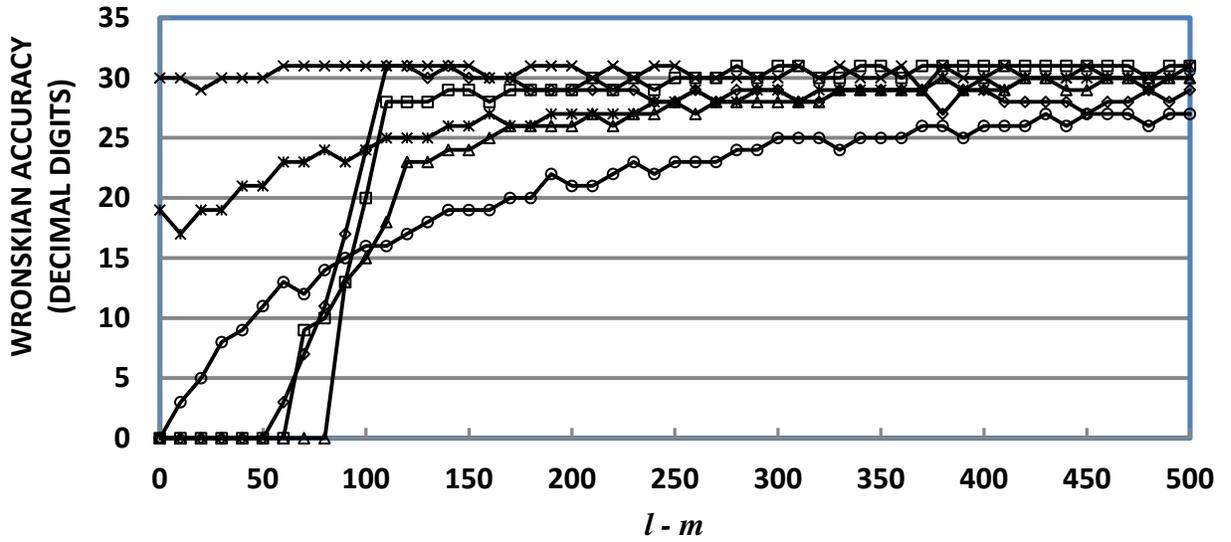

Fig. 5: Wronskian accuracy in decimal digits obtained using the traditional Legendre function expression to calculate $R^{(2)}_{ml}(-ic, i\xi)$ for selected parameters ($\xi$, $c_r$, $c_i$, $m$): × (0.1, 100, 0, 100); * (0.1, 100, 20, 100); o (0.1, 100, 40, 100). ◊ (0.99, 100, 0, 50); □ (0.99, 100, 20, 50); Δ (0.99, 100, 40, 50)

## 10  Alternative Legendre function expression

Baber and Hasse [11] provided the following expression for the oblate radial functions of the third kind $R^{(3)}_{ml} = R^{(1)}_{ml} + iR^{(2)}_{ml}$ in terms of the functions $Q^m_{n+m}(i\xi)$:

$$R^{(3)}_{ml}(-ic, i\xi) = \frac{e^{ic\xi} i^{2m-l}}{m!c} \sum_{n=-l}^{\infty} \frac{A^{ml}_n}{A^{ml}_{-m}} Q^m_{m+n}(i\xi). \tag{17}$$

.

The coefficients $A^{ml}_n$ satisfy the following recursion relation:

$$\frac{2c(n+m+1)(n+2m+1)}{(2n+2m+3)}A_{n+1}^{ml} - [(n+m)(n+m+1) - \lambda_{ml} - c^2]A_n^{ml}$$
$$-\frac{2cn(n+m)}{(2n+2m-1)}A_{n-1}^{ml} = 0,$$
(18)

with the asymptotic condition $A_{n+1}^{ml}/A_n^{ml} \to c/n$ as $n \to \infty$. Flammer [2, p. 40] and [12] provide a discussion of this expression. The radial functions of the second kind are then given by the imaginary part of the right hand side of (17). When $c$ is real, this expression does not suffer subtraction errors for the lowest order functions when both $c$ and $m$ are small to moderate in size and $\xi$ is not large. It can provide accurate results at lower values of $l - m$ when the traditional Legendre function expression fails to do so. However, there are new alternative methods available in coblfcn that appear to do this. The alternative Legendre function method will be included in coblfcn but rarely if ever used.

## 11  Integral expressions for calculating radial functions

Reference [10] shows that the integral expressions given by Flammer [2, pp.53-54] are useful for calculating the prolate radial functions of the second kind. Converted to oblate form, these expressions become:

$$R_{ml}^{(2)}(-ic,i\xi) = \frac{(-1)^{(l-m)/2}(2m+1)}{2^{m+1}m!d_0(-ic|ml)} \times$$
$$\int_{-1}^{+1}\left[\frac{(\xi^2+1)(1-\eta^2)}{(\xi^2-\eta^2+1)}\right]^{m/2} y_m[c(\xi^2-\eta^2+1)^{1/2}]S_{ml}^{(1)}(-ic,\eta)d\eta, \quad l-m \text{ even.}$$
(19)

$$R_{ml}^{(2)}(-ic,i\xi) = \frac{(-1)^{(l-m-1)/2}(2m+3)}{2^{m+1}m!d_1(-ic|ml)} \times$$
$$\int_{-1}^{+1}\frac{[(\xi^2+1)(1-\eta^2)]^{m/2}}{(\xi^2-\eta^2+1)^{(m+1)/2}}\xi\eta\, y_{m+1}[c(\xi^2-\eta^2+1)^{1/2}]S_{ml}^{(1)}(-ic,\eta)d\eta, \quad l-m \text{ odd.}$$
(20)

It is convenient to define $z = c(\xi^2 - \eta^2 + 1)^{1/2}$ and a window function $F_m(\xi,\eta)$ given by

$$F_m(\xi,\eta) = \left[\frac{(\xi^2+1)(1-\eta^2)}{\xi^2-\eta^2+1}\right]^{m/2}.$$
(21)

Expanding $S_{ml}^{(1)}$ in (19) and (20) in terms of associated Legendre functions results in:

$$R_{ml}^{(2)}(-ic,i\xi) = B_{ml}^{(a)}(c)\sum_{n=0}^{\infty}{}'d_n(-ic|ml)I_{mn}^{(a)}(c,\xi), \quad l-m \text{ even,} \tag{22}$$

$$R_{ml}^{(2)}(-ic,i\xi) = \xi B_{ml}^{(b)}(c)\sum_{n=1}^{\infty}{}'d_n(-ic|ml)I_{mn}^{(b)}(c,\xi), \quad l-m \text{ odd,} \tag{23}$$

where $B_{ml}^{(a)}(c)$ is the leading coefficient in (19), $B_{ml}^{(b)}(c)$ is the leading coefficient in (20), and

$$I_{mn}^{(a)}(c,\xi) = \int_{-1}^{+1}F_m(\xi,\eta)y_m(z)P_{m+n}^m(\eta)d\eta, \quad l-m \text{ even,} \tag{24}$$

$$I_{mn}^{(b)}(c,\xi) = c\int_{-1}^{+1}[F_m(\xi,\eta)/z]\eta\, y_{m+1}(z)P_{m+n}^m(\eta)d\eta, \quad l-m \text{ odd.} \tag{25}$$

One obtains corresponding expressions for the first derivatives of $R_{ml}^{(2)}$ with respect to $\xi$ from (19) and (20) by differentiating, utilizing standard recursion relations for the spherical Neumann functions, and collecting terms. This gives:

$$\frac{dR_{ml}^{(2)}}{d\xi}(-ic,i\xi) = \frac{m\xi}{\xi^2+1}R_{ml}^{(2)}(-ic,i\xi) - \frac{(-1)^{(l-m)/2}(2m+1)c\xi}{2^{m+1}m!d_0(-ic|ml)} \times$$
$$\int_{-1}^{+1}\frac{[(\xi^2+1)(1-\eta^2)]^{m/2}}{(\xi^2-\eta^2+1)^{(m+1)/2}}y_{m+1}[c(\xi^2-\eta^2+1)^{1/2}]S_{ml}^{(1)}(-ic,\eta)d\eta, \quad l-m \text{ even}, \quad (26)$$

$$\frac{dR_{ml}^{(2)}}{d\xi}(-ic,i\xi) = \frac{(m+1)\xi^2+1}{\xi(\xi^2+1)}R_{ml}^{(2)}(c,\xi) - \frac{(-1)^{(l-m-1)/2}(2m+3)c\xi^2}{2^{m+1}m!d_1(-ic|ml)} \times$$
$$\int_{-1}^{+1}\frac{[(\xi^2+1)(1-\eta^2)]^{m/2}}{(\xi^2-\eta^2+1)^{(m+2)/2}}\eta\, y_{m+2}[c(\xi^2-\eta^2+1)^{1/2}]S_{ml}^{(1)}(-ic,\eta)d\eta, \quad l-m \text{ odd}. \quad (27)$$

Replacing $S_{ml}^{(1)}$ with its expansion in (1) results in:

$$\frac{dR_{ml}^{(2)}}{d\xi}(-ic,i\xi) = \frac{m\xi}{\xi^2+1}R_{ml}^{(2)}(c,\xi) + c\xi B_{ml}^{(a)}\sum_{n=0}^{\infty}{'}d_n(-ic|ml)I_{mn}^{(c)}(c,\xi), \quad l-m \text{ even}, \quad (28)$$

$$\frac{dR_{ml}^{(2)}}{d\xi}(-ic,i\xi) = \frac{(m+1)\xi^2+1}{\xi(\xi^2+1)}R_{ml}^{(2)}(-ic,i\xi) + c\xi^2 B_{ml}^{(b)}\sum_{n=1}^{\infty}{'}d_n(-ic|ml)I_{mn}^{(d)}(c,\xi), \quad l-m \text{ odd}, \quad (29)$$

where

$$I_{mn}^{(c)}(c,\xi) = \int_{-1}^{+1}\left[F_m(\xi,\eta)/z\right]y_{m+1}[z]P_{m+n}^m(\eta)d\eta, \quad l-m \text{ even}, \quad (30)$$

$$I_{mn}^{(d)}(c,\xi) = \int_{-1}^{+1}\left[F_m(\xi,\eta)/z^2\right]\eta y_{m+2}[z]P_{m+n}^m(\eta)d\eta, \quad l-m \text{ odd}. \quad (31)$$

The required integrals $I_{mn}(c,\xi)$ have an integrand that is symmetric about $\eta = 0$. They can be computed using Gauss quadrature over positive values of $\eta$ and doubling the result. However, one must be careful to increase the density of quadrature points near $\eta = 1$ when $\xi$ approaches zero because of the singularity of the spherical Neumann functions at $z = 0$. The integrals tend to decrease in magnitude as $l$ increases, the decrease accompanied by loss of accuracy from increasing subtraction error. This causes a decrease in accuracy in both $R_{ml}^{(2)}$ and its first derivative as $l$ increases, although the decrease can be very gradual in many cases.

Relationships between the different integrals can be obtained through use of recursion relations for the associated Legendre functions. For example, replacing $\eta P_{m+n}^m(\eta)$ in the rhs of (31) with its equivalent in terms of $P_{m+n-1}^m(\eta)$ and $P_{m+n+1}^m(\eta)$ results in:

$$(2n+2m+1)I_{mn}^{(b)} = (n+2m)I_{m,n-1}^{(c)} + (n+1)I_{m,n+1}^{(c)}. \quad (32)$$

Thus the integrals $I_{mn}^{(b)}$ can be calculated directly from (32) instead of computing them using Gauss quadrature. Other derived relations are not as useful for calculating radial functions with a given order $m$ since they relate integrals of one kind and order $m$ to integrals of a second kind and order $m \pm 1$. These could, however, be useful when one is computing the radial functions for a range of $m$ values.

When $c$ is real, use of the integral expressions provides accurate function values over a wide range of parameters. It is especially useful when $\xi \leq 0.2$, $c$ is large, $l - m$ is less than the breakpoint but above the region where one can take advantage of eigenvalue pairing to obtain sufficiently accurate values for $R_{ml}^{(2)}$ and its first derivative. Here other methods usually suffer too much subtraction error to allow accurate results. When $\xi$ is very small and $m$ is very large, this method can often provide accurate results to very high values of $l - m$.

When $c$ is complex, the accuracy of the function values obtained using this method tends to decrease as $c_i$ increases. Figure 6 shows the Wronskian accuracy obtained when using the integral expressions for several parameter sets. Here $\xi = 0.001$, $c_r = 400$, and $m = 0$ and 100. For each $m$, there are plots of the Wronskian accuracy as a function of $l - m$ for $c_i = 0$, 20, and 40. Although the accuracy tends to decrease with increasing $c_i$, the integral method can often provide reasonably accurate values for $R_{ml}^{(2)}$ and it first derivative for complex $c$ when none of the other methods do. The spikes in the accuracy for $c_i = 20$ and 40 occur where the prolate-like eigenvalues are located. The tendency of the accuracy to be maintained to high values of $l - m$ seen in Fig. 6 results from the small value for $\xi$. As $\xi$ increases, the accuracy can fall off earlier with increasing $l - m$, especially for low values of $m$. For example, when $\xi = 0.1$, the integral expressions provide accuracy for $c = 400 + 20i$ and $m = 0$ equal to 23 digits at $l = 200$, 16 digits at $l = 300$, 10 digits at $l = 400$ and 5 digits at $l = 500$.

Use of the Wronskian tends to overestimate the accuracy of $R_{ml}^{(2)}$ for $l - m$ even and its first derivative for $l - m$ odd when using the integral expressions. Coblfcn adjusts the Wronskian estimate downward based on subtraction errors in their calculation to obtain a better estimate of accuracy. The integral expressions are not used in coblfcn when $\xi$ is greater than 0.2. Other expressions work well in that case. They are also not used for $\xi$ less than 0.0005 unless $c_i$ is greater than 2, where they can be used for values of $\xi$ as small as 0.0000001.

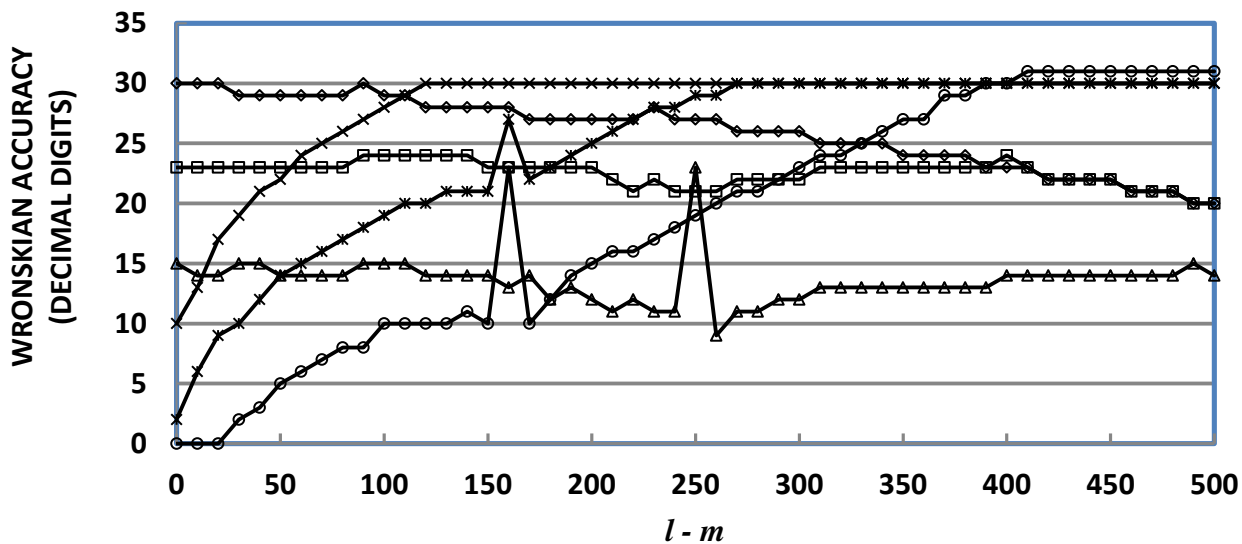

Fig. 6: Wronskian accuracy in decimal digits obtained using the integral method to calculate $R_{ml}^{(2)}(-ic, i\xi)$ for selected parameters ($\xi$, $c_r$, $c_i$, $m$): ◊ (0.001, 400, 0, 0); □ (0.001, 400, 20, 0); Δ (0.001, 400, 40, 0); × (0.001, 400, 0, 100); * (0.001, 400, 20, 100); o (0.001, 400, 40, 100).

# 12 Calculation of $R_{ml}^{(2)}(-ic, i\xi)$ using Neumann function expressions

Consider first the traditional Neumann function expression given in (9). For *c* real, this expression has the advantage that the denominator term is the corresponding angular function evaluated at $\eta = 1$, which is numerically robust with no subtraction error. However, the numerator term is asymptotic and not absolutely convergent for any finite value of $c\xi$. It can often provide accurate values, especially when $\xi$ is not small. To evaluate (9) one takes the partial sum of the series up to and including the term where the magnitude of the relative contribution is smaller than $10^{-ndec}$. The integer *ndec* is the number of decimal digits that are available in the arithmetic used in the calculations. The corresponding expression for the derivative of the radial function behaves similarly to (9). Sometimes the relative contribution never gets as small as $10^{-ndec}$. In that case the series is truncated at the term where the relative contribution is minimum. After [1] was published, it was realized that the traditional Neumann function expressions are not needed in oblfcn or in coblfcn. The subroutine for this purpose was removed from oblfcn. Other methods work just as well in regions appropriate for these expressions.

The alternative Neumann function expressions (10) and (11) obtained when $\eta$ has been set equal to 0 are very useful. The numerator sums behave as if they were not asymptotic. They are well-behaved and converge to the desired accuracy, even at high values of *c* and low values of *l - m*. There is no evidence of the series beginning to diverge as further terms are added, even when tens of thousands of additional terms are taken in the series. This is true for values of $\xi$ as low as 0.01, although the number of terms required is much larger for lower values of $\xi$. Now the denominator in (10) and (11), which is the Flammer normalization sum, does suffer the subtraction error of the corresponding angular function at $\eta = 0$. As seen in Fig. 3 in [1] for *c* real, the error is greatest at $l = m$ and decreases with increasing *l - m* until it reaches zero near the breakpoint $n_b$. For given *l - m* it increases with increasing *c* and decreases with increasing *m*. Figure 2 in [1] shows this behavior for *c* real for the case $l = m$. When *c* is complex, the behavior is nearly the same up to the breakpoint where the subtraction error is 0. Beyond this point the subtraction increases slowly with increasing *l*, reaches a secondary maximum and then decreases slowly to zero. The size of the maximum increases approximately linearly with $c_i$ and increases slowly with increasing $c_r$. It appears to be nearly independent of *m*. Figure 7 illustrates this behavior for selected parameter values.

The numerator sums in (10) and (11) also suffer subtraction errors that are similar to this. The subtraction errors restrict the use of (10) and (11) to values of *l - m* that are large enough so that the numerator and denominator achieve the desired accuracy for the radial functions. Figure 11 in [1] shows the accuracy of $R_{ml}^{(2)}$ and its first derivative obtained using (10) and (11) when *c* is real. Figure 8 below shows the corresponding accuracy for selected parameters with *c* complex. When *c* is real, the accuracies obtained using (10) and (11) are not very dependent on the value for $\xi$. This is obviously not the case for complex *c* as seen by comparing the curves for $\xi = 0.1$ and $\xi = 1.5$ when c = 200 + 40i. The dip in the curves for intermediate values of *l - m* is due to the secondary maximum that occurs in the calculation of the subtraction error for the denominator, as seen in Fig. 7. The use of this method is restricted to $\xi \geq 0.01$

When the Neumann functions $y_n(z)$ have an argument with an imaginary part that is not very small, one must take care in computing them. For real arguments, forward recursion of the standard expression starting with values for $y_0$ and $y_1$ provides accurate values for higher order functions. When the imaginary part is not small, the Neumann functions decrease in magnitude. from $y_0$ until they reach a turning point where their magnitude is approximately equal to $1/z$. See the discussion above regarding (14) and (15). Above this point they increase in magnitude with increasing order. Forward recursion from orders 0 and 1 would result in many inaccurate Neumann functions.

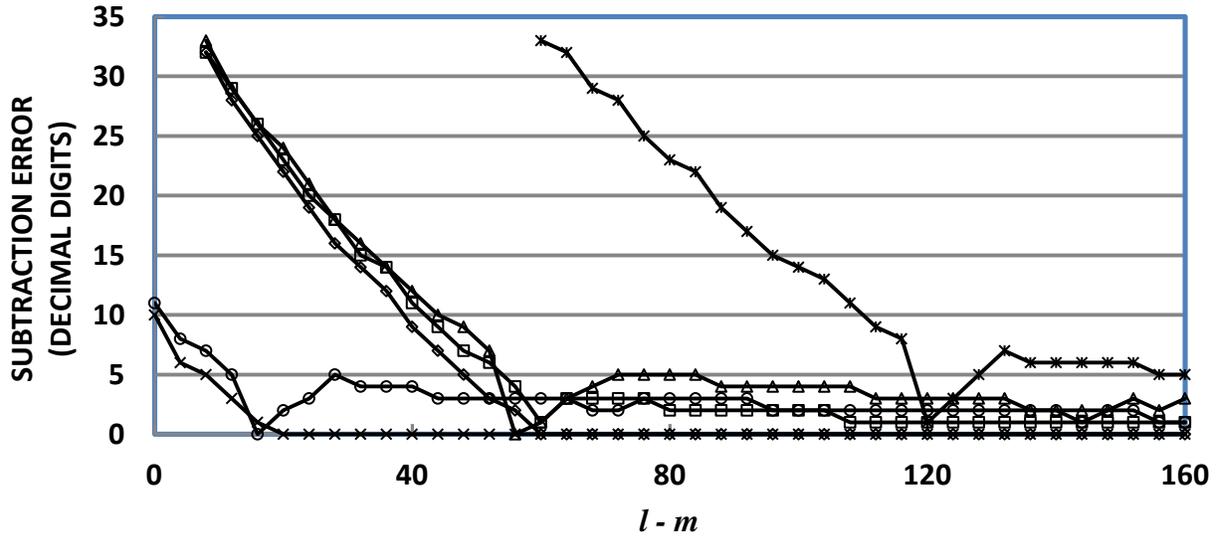

Fig. 7: Subtraction error when calculating the Flammer normalization for selected parameters ($c_r$, $c_i$, $m$): ◊ (100, 0, 0); □ (100, 20, 0); Δ (100, 40, 0); * (200, 40, 0); × (100, 0, 50); o (100, 40, 50).

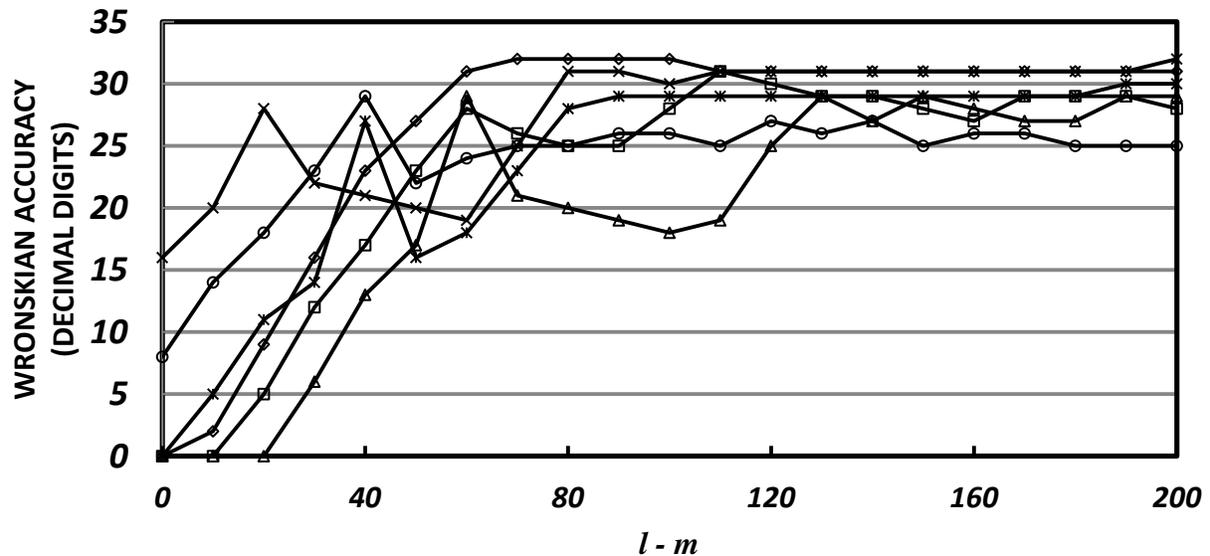

Fig. 8: Wronskian accuracy in decimal digits obtained using the Neumann function expressions with $\eta = 0$ to calculate $R_{ml}^{(2)}(-ic, i\xi)$ for selected parameters ($\xi$, $c_r$, $c_i$, $m$): ◊ (0.5, 100, 0, 0); □ (0.5, 100, 20, 0); Δ (0.5, 100, 40, 0); × (0.5, 100, 40, 50); * (0.1, 200, 40, 100); o (1.5, 200, 40, 100).

A solution for this is to first use backward recursion of the corresponding Bessel functions ratios $j_{n+1}/j_n$ from a sufficiently high value of $n$ where the ratio can be set equal to 0. The resulting ratios for orders below the turning point are fully accurate or nearly so. Individual Bessel function values are then obtained by forward multiplication of the ratios starting with the value for $j_0$. One then obtains values for $y_n$ below the turning point from expressions relating them to the corresponding Bessel functions. Then the standard recursion can be used to compute the Neumann functions above the turning point. For convenience, coblfcn computes ratios of successive Bessel and Neumann functions instead of the individual functions. Note that the same care described here must be taken when using the integral method with its Neumann functions kernels.

When $c_r$ is very large and $l$ is below the breakpoint, sometimes the methods described above are unable to provide sufficiently accurate values for the radial functions of the second kind and their first derivatives. This is especially true when $m$ is neither small nor extremely large. The use of a variable $\eta$ method similar to that described above in Sec. 6 for calculating $R_{ml}^{(1)}$ can often help here when $\xi$ is greater than about 0.05. It can bridge the gap in $l$ - $m$ where eigenvalue pairing is no longer sufficient, the integral method fails to provide the desired accuracy and the Neumann expression with $\eta = 0$ has not yet begun to proved accurate results.

The variable $\eta$ method for calculating $R_{ml}^{(2)}$ and its first derivative proceeds as follows. At the lowest value of $l$ - $m$ where the desired accuracy is not achieved using other methods, the value of $\eta$ is decreased in steps from unity and the radial functions of the second kind are calculated at each step. The accuracy initially tends to increase with decreasing $\eta$ as the numerator series becomes more accurate, although it may take several steps before any significant increase is obtained. Ideally, the desired accuracy is achieved after one step. The associated value of $\eta$ for that step is then used for the next value of $l$ - $m$. It continues to be used for progressively higher values of $l$ - $m$ until the accuracy again falls below the desired minimum. Then $\eta$ is again decreased in steps until the desired accuracy is achieved for that value of $l$ - $m$. Typically only one or two steps are needed here. The process is repeated until the $\eta = 0$ expression offers sufficient accuracy. Sometimes no value of $\eta$ provides the desired accuracy, even when the process is continued until the denominator becomes less accurate than the numerator series and the accuracy starts to decrease. The best $\eta$ for this value of $l$ - $m$ is the one used for the previous step, which is then used for the next value of $l$ - $m$. The process is continued until the $\eta = 0$ expressions offers the desired accuracy and is used instead. When $\xi \leq 0.99$, the traditional Legendre expressions will often replace the $\eta = 0$ expressions when $l$ - $m$ is large enough so that they provide the desired accuracy. As in Sec. 6, it is convenient to use steps in $\theta = \arccos(\eta)$. A step size of about 0.1 radian is used when $\xi > 0.4$ and about 0.05 radian for lower values of $\xi$.

Figure 9 shows some examples of the accuracy achieved using this method plotted versus $l$ - $m$. As before, the calculations were carried out with 33 decimal digits of precision. The first three curves are for $c_r = 300$, $\xi = 0.5$, and $m = 0$. The first one is for $c$ real and the second and third ones are for $c_i = 40$ and 80. Here it is seen that the accuracy for very low $l$ - $m$ is actually higher for $c_i = 40$ than for real $c$ and even higher for $c_i = 80$. However, as $l$ - $m$ continues to increase, the accuracy for $c$ real quickly reaches and maintains 30 or so digits. And the accuracy for $c_i = 40$ and 80 decreases considerably to a minimum at intermediate values of $l$ - $m$ before it eventually increases to 25 or more digits. The size of the decrease for $c_i = 80$ is about twice that for $c_i = 40$. Curves 4 and 5 are for $c_r = 300$, $\xi = 1.0$, and $m = 100$ with $c_i = 40$ for curve 4 and 80 for curve 5. Here we see that the decrease in accuracy is somewhat less than for $\xi = 0.5$. Calculations show that the minimum does not change much with $m$. Also, it tends to be located near the same value for $l$ as $m$ varies. The final curve shows the effect of reducing $\xi$ to 0.3. Here the decrease is even greater than that seen for $\xi = 0.5$. The use of the variable $\eta$ method is restricted to values of $\xi \geq 0.05$.

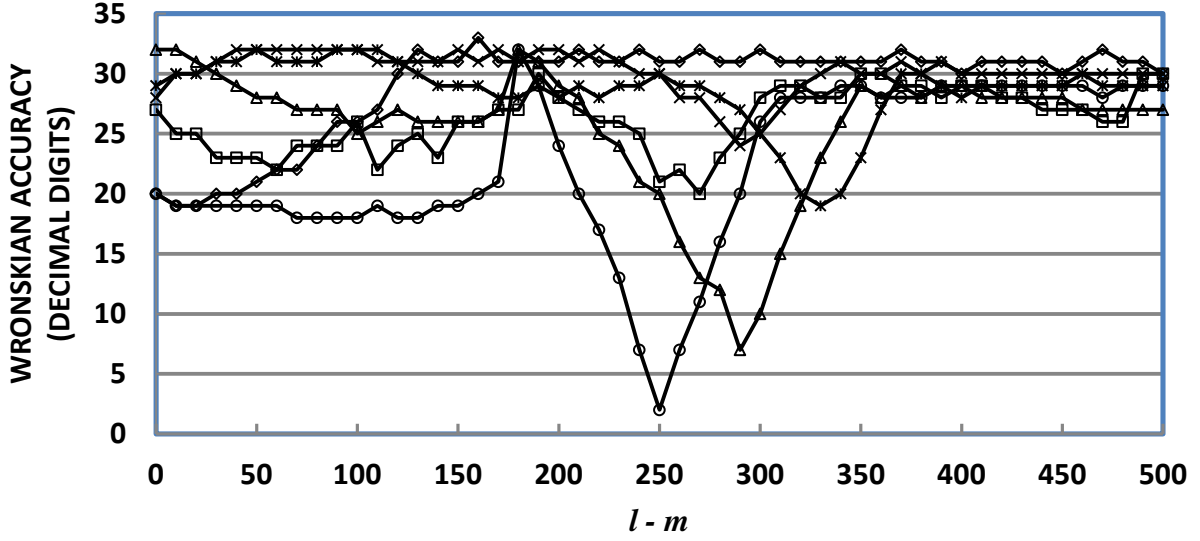

Fig. 9: Wronskian accuracy in decimal digits obtained using the variable $\eta$ method for selected parameters ($\xi$, $c_i$, $c_r$, $m$): ◊ (0.5, 300, 0, 0); □ (0.5, 300, 40, 0); Δ (0.5, 300, 80, 0); × (1.0, 300, 40, 100); * (1.0, 300, 80, 100); o (0.3, 300, 80, 0).

## 13 Calculation of radial functions for $\xi = 0$

When $\xi = 0$ and $c$ is real, accurate values for $R_{ml}^{(1)}$ are given when $l - m$ is even by the $d_0$ term in (9) since the remaining terms vanish. It is equal to zero for $l - m$ odd since all of the terms vanish. The first derivative of $R_{ml}^{(1)}$ for $l - m$ odd is given by the $d_1$ term in the derivative of (9) since the remaining terms vanish. It is equal to zero for $l - m$ even since all of the terms vanish. The Wronskian can be used to obtain accurate values for $R_{ml}^{(2)}$ for $l - m$ odd and for its first derivative when $l - m$ is even. Here

$$R_{ml}^{(2)}(-ic,0) = -\frac{1}{c\dfrac{dR_{ml}^{(1)}(-ic,0)}{d\xi}}, \qquad l-m \ odd,$$

$$\frac{dR_{ml}^{(2)}(-ic,0)}{d\xi} = \frac{1}{cR_{ml}^{(1)}(-ic,0)}, \qquad l-m \ even.$$

(33)

When $c$ is complex, subtraction errors in computing the Morse-Feshbach normalization will reduce the accuracy of the nonzero radial functions of the first kind and the radial functions of the second obtained using (33). These subtraction errors have a maximum value where the lowest order prolate-like occur. This is typically somewhat below the breakpoint. Figure 2 shows this maximum value for a range of $c_i$. The accuracy of $R_{ml}^{(1)}$, its first derivative and $R_{ml}^{(2)}$ or its first derivative obtained using (33) is conservatively estimated using the subtraction error in the Morse and Feshbach normalization, the estimated accuracy of the eigenvalue, and the number of digits of match between the forward and backward recursions for the $d$ coefficients.

Values for $R_{ml}^{(2)}$ when $l - m$ is even are obtained from the limiting form of (16) while values for its first derivative for $l - m$ odd are given by the limiting form of the first derivative of (16). These limiting forms suffer the subtraction errors incurred in calculating the joining factor as discussed above in Sec. 9. As such, they can be much less accurate than the radial functions of the first kind and the radial functions of the second kind obtained using (33). Figure 10 shows the subtraction error for the joining factor plotted versus $l - m$ for selected values of $c_r$, $c_i$, and $m$. The errors at low $l - m$ are similar to those for real $c$. The spikes in the curves correspond to the location of the prolate-like eigenvalues. When $c$ is real or $c_i$ is small and $l - m$ decreases from the breakpoint, the resulting values for $R_{ml}^{(2)}$ when $l - m$ is even and its first derivative when $l - m$ is odd become smaller in magnitude than the corresponding non-zero values for $R_{ml}^{(1)}$ or its first derivative at the same time they become less accurate. In this case, their contribution to the solution of problems involving oblate spheroidal geometry is proportionately reduced along with their reduced accuracy. The accuracy of $R_{ml}^{(2)}$ for $l - m$ even and its first derivative for $l - m$ odd is conservatively estimated using the subtraction error in the joining factor, the estimated accuracy of the eigenvalue, and the number of decimal digits of agreement between the forward and backward recursions to obtain the $d$ coefficients. When the estimated accuracy is 0, $R_{ml}^{(2)}$ for $l - m$ even and its first derivative for $l - m$ odd are set equal to zero..

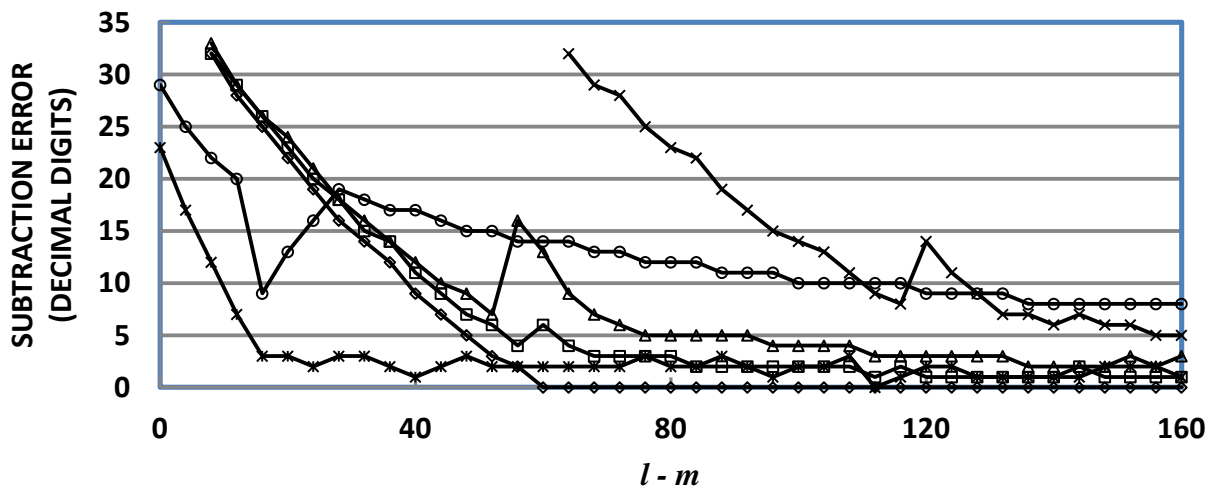

Fig. 10: Subtraction error when calculating the joining factor for selected parameters ($c_r$, $c_i$, $m$): ◊ (100, 0, 0); □ (100, 20, 0); Δ (100, 40, 0); * (100, 0, 50); o (100, 40, 50); × (200, 40, 0).

## 14 A Fortran program to compute oblate spheroidal functions for complex argument

A Fortran computer program called coblfcn [13] has been developed to calculate the oblate spheroidal functions when the size parameter *c* is complex. It can also calculate the functions if *c* is real, but an earlier program oblfcn for this purpose runs faster. Coblfcn is available as either a stand-alone program or as a subroutine. In the stand-alone program, the calculated radial and/or angular functions are written to the files fort.20 and/or fort.30, respectively. In the subroutine version, the function values are passed via the subroutine statement. Coblfcn performs calculations in either double precision arithmetic, quadruple precision arithmetic or in a hybrid mode where double precision is used for all but the Bouwkamp procedure, which is performed in quadruple precision. The choice of arithmetic is controlled by a module called param at the end of coblfcn that sets the kind parameters knd and knd1 equal to the number of bytes for real data in either double precision or quadruple precision. The parameter kind1 applies just to the Bouwkamp procedure. Coblfcn was developed on a laptop where kind equaled either 8 or 16. It provided an approximate precision of either 15 or 33 decimal digits. Calculation options include (1) radial functions of the first kind and their first derivatives, (2) radial functions of both the first and second kind and their first derivatives, (3) angular functions of the first kind, and (4) angular functions of the first kind and their first derivatives. If desired, both radial and angular functions can be calculated during the same run.

As discussed above in Sec. 2, coblfcn provides an estimate of the number of accurate digits in the angular functions and their first derivatives based on the subtraction errors involved in their calculation and normalization. It also provides an estimate of the number of accurate digits in $R_{ml}^{(2)}$ and its first derivative. Both $R_{ml}^{(1)}$ and its first derivative are almost always highly accurate. The estimate of accuracy is often based on the Wronskian. When $\xi = 0$, the estimate of accuracy is based on either the subtraction error in either the Morse and Feshbach normalization or the joining factor. This is discussed in Sec. 13. The output of coblfcn includes diagnostic files (fort.40 for radial functions and fort.30 for angular functions) including information such as the number of terms both available and used in the various series. The generation of these files can be suppressed if desired. The default mode for the subroutine version is for suppressing them. Coblfcn also provides a file fort.60 that tells the user when the estimated accuracy for the radial functions falls below a specified number of decimal digits. It alerts the user when the Bouwkamp procedure produces the same eigenvalue for two values of *l* with the same parity. Here the results from coblfcn for the value of *m* where it occurs are not useful. This should only happen when coblfcn is used outside its usual parameter ranges. See below for a discussion of these ranges. Fort.60 also alerts the user when the estimated accuracy of the Meixner and Schäfke normalization is less than the same specified number of digits. Here the values of *m, l*, and *c* where this occurs are written to fort.60.

Coblfcn calculates angular functions of the first kind using (1) and gives them unit norm or normalizes them using the Meixner and Schäfke normalization, depending on the input parameter iopnorm. It does so for a specified number of *l* values beginning with $l = m$. In the stand-alone program it allows either $\eta$ or $\theta = \arccos(\eta)$ arguments and computes angular functions for a range of arguments determined by a first value, an increment and the number of arguments desired. It does so for a range of *m* values determined by a first value, an increment and the number of *m* values desired. The resulting function values are given as a characteristic with a magnitude between 1.0 and 10.0 and an integer exponent iexp that denotes the power of 10 for the factor $10.0^{iexp}$. In the subroutine version, angular functions values are obtained for a single value of *m* and an input vector of $\eta$ values. It is expected that the user will choose unit norm to avoid any potential overflow problems for very high values of *m* if the characteristics and exponents are combined.

Coblfcn calculates radial functions for a single input value of $\xi$ and for a specified number of values for the degree *l* beginning with *l* = *m*. In the stand-alone program it does so for a range of orders *m* values, while in the subroutine version it does so only for a single value of *m*. Coblfcn calculates the radial functions of the first kind using the traditional Bessel function expression together with the variable $\eta$ method described above. The radial functions of the first kind are almost always highly accurate (unless near a root). Coblfcn obtains radial functions of the second kind using either the pairing of low-order eigenvalues via (12), the near equality of $i\,R_{ml}^{(1)}$ and $R_{ml}^{(2)}$ for large values of $c_i\xi$ as discussed in Sec. 8, the alternative $\eta = 0$ Neumann function expressions (10) and (11), the traditional associated Legendre function expression (16), the Baber and Hasse Legendre function expression (17), the variable $\eta$ method, or the integral expressions in Sec. 11. The methods used in coblfcn are based on the input parameters and the integer minacc that specifies the number of accurate decimal digits that are desired. Minacc is set equal to 8 for double precision arithmetic. For quadruple precision it is set equal to 15 unless $c_i$ is greater than 20, where it is set equal to 8. The value for quadruple precision arithmetic can be changed if desired, especially if higher accuracy is desired for input parameters where it can be achieved. It is advised to leave minacc set to 8 for double precision arithmetic unless the precision available in double precision on the user's computer is greater than 15.

The methods used in coblfcn to calculate the radial functions of the second kind are based on the input parameters and the desired minimum accuracy. Coblfcn starts at *l* = *m* with the use of paired eigenvalues if the pairing is sufficient to provide the desired accuracy and continues with increasing *l* - *m* until the pairing is insufficient. If $c_r$ is not very large, then coblfcn tends to use the alternative $\eta = 0$ Neumann function expressions for $\xi$ > 0.99 and the traditional Legendre function expression for $\xi \leq 0.99$.

For larger values of *c* and smaller values of $\xi$, coblfcn switches to the integral method after the paired eigenvalue method and continues until the traditional Legendre function expansion provides sufficient accuracy and is used for all higher values of *l* - *m* that are desired. For larger values of both *c* and $\xi$, it switches to one of the Neumann function expansions and continues until the $\eta$ = 0 Neumann expansion provides the desired accuracy and is used for all higher values of *l* - *m* that are desired. When $c_i$ is large and $R_{ml}^{(1)}$ is large in magnitude, the near equality of $i\,R_{ml}^{(1)}$ and $R_{ml}^{(2)}$ is used when it provides the minimum desired accuracy or if it provides more accuracy than the other methods.

The various series involved in coblfcn are computed in a way to avoid potential overflow and underflow in the calculations. First the expansion coefficients and the expansion functions are calculated as ratios using appropriate recursion relations. See comment statements in the appropriate subroutine in coblfcn where they are calculated. Bessel and Neumann function series are summed starting with the *n* = *l* − *m* term while Legendre function and integral method series are summed starting with the lowest term. The function and coefficient values for the first term are factored out of the expansion and the first term is set equal to unity. Summation is performed using ratios to obtain the next terms in the series. For Bessel and Neumann functions expansions summations are taken both forward and backward. The resulting value is then multiplied by the first term to obtain the desired sum. Here the relevant Bessel, Neumann, or Legendre function for the first term taken has been computed previously by forward multiplication of ratios starting with known values for the lowest two values. For example, this would be $j_0$ and $j_1$ for Bessel functions. During forward multiplication, power 10 exponents are stripped out of the product at each step to avoid either underflow or underflow. This results in a characteristic with magnitude between 1.0 and 10.0 and an exponent denoting the corresponding power of 10. Radial function values are stored as both a characteristic and an exponent. This allows coblfcn to provide results at high values of *l - m* where the radial functions of the first kind would otherwise underflow and the radial functions of the second kind would overflow.

## 15    Estimated accuracy

Coblfcn was tested extensively using a laptop pc and a Fortran compiler that provides approximately 15 decimal digits in double precision (64 bit) arithmetic and approximately 33 digits in quadruple precision (128 bit) arithmetic. If the user's computer provides a different number of digits, the following estimates should be adjusted up or down depending on whether more or fewer decimal digits are provided. Testing included values of $\xi$ ranging from 0.000001 to 10 as well as the special case $\xi = 0$, values for $c_r$ up to 5000, and values of $c_i$ up to 200. Testing for both the double precision and the hybrid versions included all values of the degree *m* from 0 to 200 and from 210 to 1000 in steps of 10. Testing for the quadruple precision version included values of *m* from 0 to 200 in steps of 10 and from 250 to 1000 in steps of 50. For all three versions, the values of the degree *l* ranged from *m* to a value high enough so that the magnitudes of $R_{ml}^{(1)}$ and its first derivative were less than $10^{-300}$.

The Wronskian (13) usually provides the estimated accuracy of the radial functions. There are several exceptions. First is when $c_i\xi$ is large and $R_{ml}^{(1)}$ is sufficiently large so that subtraction errors occur in forming the Wronskian, as described in Sec. 8. Second is when the Wronskian is used to improve the accuracy of the results from the Legendre function expressions, as described in Sec. 9. The third is when the Wronskian estimate for the integral method is sometimes adjusted downward for very small values of $\xi$, as described in Sec. 11. The fourth is whenever the traditional Legendre expression is used. Comparison of double and quadruple precision results showed that for larger values of $c_i$ the Wronskian sometimes overestimated the radial function accuracy at values of *l - m* below the breakpoint using this method. Here the accuracy is given as the lesser of the Wronskian estimate and one based on subtraction errors in the series calculation of $R_{ml}^{(2)}$ and its first derivative as well as $R_{ml}^{(1)}$ and its first derivative. Finally is the case where $\xi = 0$. Here the accuracy is estimated using methods described in Sec. 13.

In the following discussion, the term useful results means that the estimated accuracy for the radial functions observed during testing never fell below 5 decimal digits unless otherwise stated. It is expected that there are many applications where occasional 5 digit results are acceptable. Possibly even an isolated 4 digit result is acceptable. Note that there is no guarantee that the estimated accuracy for parameter values other than those tested will be as high as described. The discussion below will focus on $\xi$ unequal to zero. It is expected that function values for $\xi = 0$ will be useful for the same values of *c* and *m* that useful results are obtained for at least one value of $\xi$. Here both $R_{ml}^{(2)}$ and its first derivative will be highly accurate. Whenever values of $R_{ml}^{(2)}$ for *l - m* even are less accurate than 5 digits, they are expected to be proportionally smaller in magnitude than $R_{ml}^{(1)}$. Similarly for the first derivatives of $R_{ml}^{(2)}$ and $R_{ml}^{(1)}$ when $l - m$ is odd.

Using double precision arithmetic, including for the Bouwkamp procedure, coblfcn provides useful results for $c_i$ up to 10, for $c_r$ up to 5000, for $m$ up to at least 1000 and for all tested values of $\xi$ down to 0.000001. When $c_i$ is less than about 5, the resulting accuracy is similar to that provided by oblfcn for $c$ real. Extensive testing for $c_i = 10$ showed that $R_{ml}^{(1)}$ and its first derivative are almost always accurate to 10 or more decimal digits. For all values of $\xi$ except zero, $R_{ml}^{(2)}$ and its first derivative are usually accurate to 8 or more decimal digits, but accuracies lower than this were sometimes seen, especially for larger $c_r$ and small $\xi$. Nearly all accuracies less than 8 digits occurred near but somewhat below the so-called breakpoint. No 5 digit results were seen for $c_r$ up to 200 or so. Only a few 5 digit results were seen for $\xi \geq 0.001$, even for $c_r = 5000$. The largest number of 5 digit results occurred for $c_r = 5000$ and $\xi = 0.000001$. Even here, however, there were no more than about 3 such results for each value of $m$.

Similar testing for $ci = 12$ showed a few more 5 digit results for small $\xi$. At least 5 digits were obtained for $c_r$ up to 5000 for all values of $\xi$ down to 0.000001.

Testing for $c_i = 15$ showed yet more 5 digits results and even 4 digit results for $c_r = 5000$. However, at least 5 digits were obtained for $c_r$ up to 2000 for all values of $\xi$ down to 0.000001.

Testing for $c_i = 20$ showed 5 or more digits of accuracy for $\xi$ down to 0.000001 when $c_r \leq 100$. There were 5 or more digits of accuracy at $\xi \geq 0.1$ for $c_r = 150$ and at $\xi \geq 0.2$ for $c_r$ up to 2000. For $c_r = 5000$, duplicated eigenvalues appeared for some values of $m$. This occurs because of the decreasing accuracy of the matrix estimates of the eigenvalues near the breakpoint as $c_i$ increases, especially for large values of $c_r$. Even if one uses quadruple precision for the Bouwkamp procedure to obtain convergence to an eigenvalue, the converged eigenvalue may not be the correct one, especially if the estimate is closer to another eigenvalue of the same parity than to the desired one.

Testing for $c_i = 25$ showed 5 or more digits of accuracy for $\xi \geq 0.2$ when $c_r \leq 100$, for $\xi \geq 0.3$ when $c_r = 150$, and for $\xi \geq 0.4$ for $c_r = 200$ and 250. The ranges here can be extended somewhat to higher values of $c_r$ by using quadruple precision for the Bouwkamp procedure.

Testing for yet higher values of $c_i$ showed a continued increase in the minimum value of $\xi$ and the maximum value of $c_r$ for which useful function values were obtained. If the user is interested in values of $\xi$ somewhat larger than 0.3 together with moderate to small values of $c_r$, then the double precision or hybrid version of coblfcn may be useful when $c_i$ is greater than 25. Otherwise, it will be necessary to use the quadruple precision version. It is recommended that the file fort.60 be used to assure that you are obtaining the accuracy you need for both the radial and angular functions and that there are no repeated eigenvalues of the same parity. Testing showed the appearance of duplicated eigenvalues for $c_i = 45$, even when $c_r$ was no greater than 200.

When the accuracy obtained using double precision is insufficient, much higher accuracy can be obtained using quadruple precision. However, coblfcn runs faster by a factor up to 50 or more in double precision than it does in quadruple precision. Running coblfcn on an ordinary laptop computer such as was used in its development can take a long time with quadruple precision when $c_r$ is extremely large.

Testing using quadruple precision for values of $c_i$ up to 40 showed that coblfcn provides useful results for $c_r$ up to at least 5000, for $m$ up to 1000 and $\xi$ down to 0.000001. Estimated accuracies for the radial functions were 8 or more digits except for an occasional 7 digit result or a rare 6 digit result that occurred primarily at lower values of $\xi$ and higher values of $c_r$. See the discussion above in Section 14 about estimated accuracy for $\xi = 0$.

Testing for $c_i = 50$ shows useful results for $c_r$ up to at least 2000. Accuracies for the radial functions were almost always at least 8 or more digits but occasional accuracies as low as 5 digits were seen near the breakpoint, primarily for $c_r \geq 1000$ and $m > 200$. When $c_r = 2000$, there were even a few 4 digit results for $m \geq 700$ and a few 3 digit results for $m \geq 800$. It is unlikely that values of $m$ this large will be required by the user.

Testing for $c_i = 60$ shows useful results for $c_r \leq 100$ although a few 4 digit results occurred at $\xi = 0.00001$ and $0.000001$. Results for $c_r = 200$ were similar except that there were a few more 4 digit results at $\xi = 0.00001$ and $0.00001$ and some 4 digit results now at $\xi = 0.0001$ and $0.001$. Testing for $c_r = 500$ shows that useful results are only obtained for $\xi \geq 0.05$ for all $m$ and for $\xi < 0.05$ for $m$ up to about 160. Testing for $c_r = 1000$ and 2000 showed useful results for $\xi \geq 0.2$ for all $m$ with a possible rare 4 digit result.

Testing for $c_i = 70$ shows useful results for $c_r \leq 100$ when $\xi \geq 0.01$, for $c_r = 150$ when $\xi \geq 0.05$, for $c_r = 200$ when $\xi \geq 0.1$, for $c_r = 500$ when $\xi \geq 0.3$ and for $c_r = 1000$ when $\xi \geq 0.5$.

Testing for $c_i = 80$ shows useful results for $c_r \leq 20$ when $\xi \geq 0.01$, for $c_r = 50$ when $\xi \geq 0.02$, for $c_r = 100$ when $\xi \geq 0.2$,

Testing for yet higher values of $c_i$ showed a continued increase in the minimum value of $\xi$ and the maximum value of $c_r$ for which useful function values were obtained. Again it is recommended that the file fort.60 be used to assure that you are obtaining the accuracy you need and that there are no repeated eigenvalues of the same parity.

Care must be taken in comparing radial functions values obtained using double precision with those using quadruple precision. One must look for the same eigenvalue to do the comparison. When the eigenvalues were ordered, the prolate-like eigenvalues were placed after eigenvalues either with a negative real part or paired eigenvalues with positive real part if they occurred. The accuracy of the matrix estimates used in the identification of paired eigenvalues can be so poor for double precision that one or more paired eigenvalues are not recognized. Therefore, the prolate-like eigenvalues are placed earlier in the sequence using double precision. The paired eigenvalues that were not recognized are then located after the prolate-like eigenvalues. It is interesting to note that this can lead to the radial functions of these non-recognized paired eigenvalues.to have opposite signs for double precision relative and quadruple precision. This occurs because of the limit that the radial function have as $c\xi \to \infty$ [2, p. 32]. This has no effect on the solution of problems using these functions.

Testing of coblfcn indicates that when all of the radial functions for a given value of $m$ have at least 5 digits of accuracy, the angular functions will usually have at least 5 accurate digits, except when they suffer large subtraction errors for lower values of $l - m$ and $\eta$ near 0 during calculation using (1). See Sec. 2 for a discussion of these subtraction errors. When they occur, the resulting angular functions and their first derivatives are reduced in magnitude by an amount corresponding to the subtraction error. Their magnitude in this case is corresponding smaller than angular functions for higher values of $l - m$ and/or $\eta$ not near zero. The loss in accuracy due to these subtraction errors will not likely affect numerical results for physical problems using these functions. The angular functions for higher values of $l$ or for $\eta$ near unity will have at least 5 digits of accuracy whenever the radial functions have at least 5 digits of accuracy.

The accuracy of the eigenvalue and the degree to which the forward and backward recursions match during calculation of the $d$ coefficients also affects the accuracy of the angular functions. However, a reduction in the accuracy of either of these has an equal impact on the accuracy of the angular and radial functions and will not result in the angular functions having a lower accuracy than the radial functions.

The accuracy of the Meixner and Schäfke normalization affects only the angular functions. It sets an upper bound on their accuracy. As seen in Fig. 1, subtraction errors in the Meixner and Schäfke normalization are zero for small $c_i$ and can become as large as 6 digits for $c_i = 20$ and 12 digits for $c_i = 50$ as $c_r$ increases to 5000. This loss in accuracy is not likely a problem when using double precision arithmetic with 15 decimal digits since as $c_i$ becomes larger than 20, the values of $c_r$ for which the radial functions are accurate to only 5 digits are progressively smaller, being 100 for $c_i = 25$ and less than this for higher $c_i$. Using double precision, the Meixner and Schäfke normalization should be accurate to at least 5 digits wherever the radial functions are also accurate to at least 5 digits. Nonetheless, it is recommended that the user utilize the file fort.60 to alert when the estimated accuracy of this normalization falls below a specified minimum value.

For higher values of $c_i$ when using quadruple precision, the loss of accuracy in the normalization factor is even greater. For $c_i = 60$, the loss of accuracy can be as large as 24 digits for $c_r = 2000$ and 25 digits for $c_r = 5000$. For $c_i = 70$, it is 26 digits for $c_r = 1000$. And for $c_i = 80$ it is 25 digits for $c_r = 400$. This should not be a problem using 33 decimal digits since it still allows for accuracies of at least 5 digits for the angular functions everywhere the radial functions also have an accuracy of 5 or more digits.

## 16 Summary

Procedures to calculate the oblate spheroidal functions for complex values of the size parameter *c* are provided. Most of the procedures are based on those described in an earlier paper [1] addressing calculation of the functions for real *c*. Some of the procedures addressed here are traditional ones, but many are alternative procedures necessary to provide function values over wide parameter ranges. The numerical behavior of each of the procedures together with examples is discussed. A Fortran computer program coblfcn is described that incorporates these procedures to provide useful values for the angular functions of the first kind and the radial functions of both kinds together with their first derivatives. Coblfcn can be run in either double precision, quadruple precision or a hybrid mode where double precision is used for all but the Bouwkamp procedure to obtain accurate eigenvalues. The paper presents an summary of the estimated accuracy using coblfcn in double and quadruple precision. It is seen that useful function values can be obtained over extremely wide parameter ranges as long as the imaginary part of *c* is less than about 20 for double precision and less than about 50 for quadruple precision. Useful results can be obtained for higher values of the imaginary part of *c* but the parameter ranges narrow. A listing of oblfcn together with sample output is freely available in text format on the web site listed in [13]. Both a stand-alone program and a subroutine version are provided.